\documentclass[10pt,twocolumn]{IEEEtran}

\usepackage{url}
\usepackage{ifthen}
\usepackage{stmaryrd}
\usepackage{cite}
\usepackage[cmex10]{amsmath} 
\usepackage{amsmath}
\usepackage{subcaption}
\usepackage{dsfont}
\usepackage{amsfonts}
\usepackage{mathrsfs}
\usepackage{amsthm}
\usepackage{amssymb}
\usepackage{breqn}
\usepackage{xcolor}

\newlength{\dhatheight}

\newcommand{\Q}{\mathds{Q}}
\newcommand{\R}{\mathds{R}}
\newcommand{\Z}{\mathds{Z}}

\newcommand{\doublevecl}{z,y_1^l}
\newcommand{\doublevecn}{z,y_1^n}

\newcommand{\mkv}{-\!\!\!\!\minuso\!\!\!\!-}

\usepackage{graphicx}
\usepackage{array}

\newtheorem{theorem}{Theorem}

\newtheorem{lemma}{Lemma}

\newtheorem{proposition}{Proposition}
\newtheorem{remark}{Remark}

\begin{document}

\title{On the Exponential Approximation of {Type II} Error Probability 
of Distributed  Test of Independence}

\author{Sebastian Espinosa$^*$,~\IEEEmembership{Student Member,~IEEE}, Jorge F. Silva,~\IEEEmembership{Senior Member,~IEEE}, \\ and Pablo Piantanida,~\IEEEmembership{Senior Member,~IEEE}
\thanks{S. Espinosa and J. F. Silva are with the Information and Decision System Group, Department of Electrical Engineering,  Universidad de Chile, Santiago, Chile (email: josilva@ing.uchile.cl).}
\thanks{This project has received funding from the European Union’s Horizon 2020 research and innovation programme under the Marie Skłodowska-Curie grant agreement No 792464.}
\thanks{This work of S. Espinosa was funded by the National Agency for Research and Development (ANID)/ScholarshipProgram/DoctoradoNacional/2018-21180693.}
}

\maketitle

\begin{abstract}
This paper studies distributed binary test of statistical independence under communication (information bits) constraints. While testing independence is very relevant  in various applications, distributed independence test is particularly useful for event detection in sensor networks where data correlation often occurs among observations of devices in the presence of a signal of interest. By focusing on the case of two devices because of their tractability, we begin by investigating conditions on \textsc{Type I} error probability restrictions under which the minimum \textsc{Type II} error  admits an exponential behavior with the sample size. Then, we study the finite sample-size regime of this problem. We derive new upper and lower bounds for the gap between the minimum \textsc{Type II} error and its exponential approximation  under different setups, including restrictions imposed on the vanishing \textsc{Type I} error probability. Our theoretical results shed light on the sample-size regimes at which approximations of the \textsc{Type II} error probability via error exponents became informative enough in the sense of predicting well the actual error probability. We finally discuss an application of our results where the gap is evaluated numerically, and we show that exponential approximations are not only tractable but also a valuable proxy for the \textsc{Type II} probability of error in the finite-length regime. 
\end{abstract}


\section{Introduction}
Motivated by decision-making problems over networks, researchers within the field of statistical  signal processing  have been involved in a wide range of research initiatives studying decision and inference problems in the presence of quantization or measurement noise or data corruption by various types of perturbations. In real-world applications, these sources of degradation come from factors such as noise observations at the sensors, communication restrictions between sensors and decision agents, or by the presence of external sources of perturbations corrupting data~\cite{Mahler_2007}. The emerging field of Internet of Things (IoT) brings new dimensions and technical challenges to the  classical problem as data is no longer centrally available at the decision end. A related emerging domain is known as signal processing in the context of unlabeled or unordered data \cite{marano_2019,Wang_2018,Zhu_2017,marano_2016,unnikrishnam_2018,liu_2018,Haghighasthoar_2018}. Another important domain, which is the general focus of this work, is distributed detection under data-compression~\cite{ahlswede1986hypothesis,han1989exponential,watanabe2018neyman}. The derivation of performance limits and characterization of statistical properties of optimal detectors have been active research areas over the past years. 

Distributed detection, data fusion, and multisensor integration have a long history in statistics, signal, and information processing at large. Fundamental works can be traced back to \cite{10.5555/524893}, \cite{4102537}, and \cite{395226}, among others. Applications of the decentralized decision framework arise in communications and sensor networks, for instance, in the context of Multiple Access Channels (MAC) \cite{5722050} and wireless sensor networks \cite{4205085}. These works do not only investigate practical solutions but, importantly, they study  theoretical guarantees and performance bounds to understand the intrinsic complexity of these problems. In \cite{4099562}, the authors derived performances in the form of error exponents of the \textsc{Type I} and \textsc{Type II} error probabilities over Fading MACs. However, explicit communication restrictions between the sensors and the fusion center still remain a challenging problem~\cite{5722050,4205085}, which  implies understanding how (detection) performances are affected by the introduction of non-trivial communication restrictions. Indeed, a crucial case of particular interest is when the fusion center receives quantized descriptions of the measurements taken by remote sensors \cite{7405352,698826}. Some recent contributions have explored the asymptotic  performance limits based on error exponents of distributed scenarios with multiple decision centers and rate constraints between sensors and detectors \cite{wigger2016testing}, \cite{xiang2012interactive}, \cite{salehkalaibar2018hypothesis}. Communications constraints have also been studied within the framework of Bayesian detection  in \cite{ 4558051} and \cite{ 6203607}.

This paper investigates  the problem of distributed binary Hypothesis Testing (HT) of   statistical independence   under   communication   (information  bits)   constraints. In particular, we focus on non-asymptotic performance bounds. More specifically, we  study the gap between the minimum \textsc{Type II}  error  probability  and  its  exponential  approximation  under restrictions on the vanishing \textsc{Type I} error  probability.  To this end, we revisit the distributed scenario with communication constraints first introduced in \cite{ahlswede1986hypothesis}. This problem consists in testing against independence where the observations (e.g., sensor measurements) come from two modalities (e.g., two sensors), as shown in Fig.~\ref{hypo}. One of the modalities is to be transmitted to the  decision-maker (or detector) using an error-free communication channel that introduces a positive rate-constraint (in bits per sample).  \cite{ahlswede1986hypothesis} derives the characterization of asymptotic performance bounds in terms of a closed-form expression  for the error exponent of \textsc{Type II} error probability given a fixed restriction on the \textsc{Type I} error probability ($\epsilon>0$)~\cite[Ths. 2 and 3]{ahlswede1986hypothesis}. Notably, the results show the effect of the communication constraints in asymptotic performance (error exponent),  which is shown to be  independent of $\epsilon$.  Later on,~\cite{han1989exponential} derive an asymptotic lower bound for the error exponent when \textsc{Type I} restriction (as a sequence) tends to zero (with the sample size $n$) at an exponential rate given by $\mathcal{O}(e^{-nr})$.
 
 

\subsection{Summary of contributions}




Our work advances state-of-the-art in very different ways.

\begin{figure*}
\centering
    \includegraphics[width=0.65\linewidth]{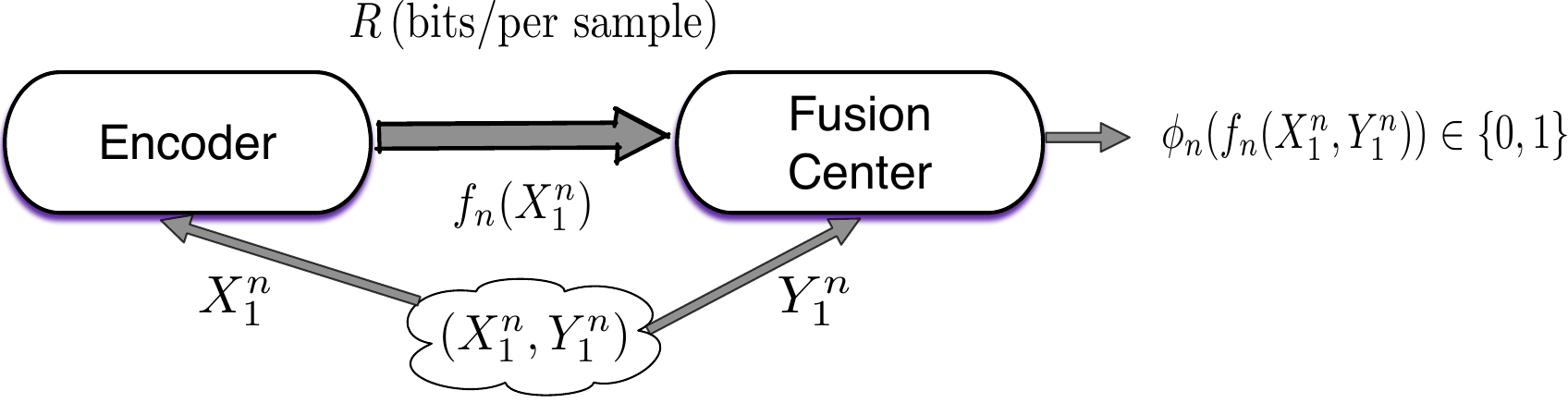} 
    \caption{Illustration of the coding-decision problem with one-side communication constraint. $f_n$ is the encoder of $X_1^n$ (one of the modalities) and $\phi_n$ is the detector acting on the one-side compressed measurements $(f_n(X^n_1),Y^n_1)$.}
      \label{hypo}  
\end{figure*}

\begin{enumerate}
    \item We study a broader family of problems (see Fig.~\ref{hypo}) where the \textsc{Type I} error probability vanishes with the sample size. The objective here is to assess the impact of this stringer set of restrictions on the asymptotic limit  of \textsc{Type II} error probability given by the error exponent. Building on concentration inequalities and results from rate-distortion theory, our main result here (cf. Theorem \ref{overheadR}) gives new conditions on the admissible converge rate of the \textsc{Type I} error probability restriction  under which the error exponent of the \textsc{Type II} error probability  admits a closed-form expression. Interestingly, for a family of sub-exponential decreasing \textsc{Type I} error probability  restrictions,  we show that the resulting error exponent matches the expression   in~\cite[Theorem 3]{ahlswede1986hypothesis} while being consistent with the results obtained for the classical communication-free problem~\cite{nakagawa1993converse}.
    
\item Regarding the  non-asymptotic regime of this problem, Theorem \ref{theorem2} offers new upper and lower bounds for the \textsc{Type II} error probabilities as a function of the number of samples, the underlying  distributions, and the restriction on the \textsc{Type I} error probability. As an important corollary, our bounds shed light on the velocity at which the error exponent is achieved as the number of samples tends to infinity, and consequently, how well the performance limits represent the performances of practical decision schemes operating based on a finite number of samples.  

\item Finally, we evaluate   our bounds numerically and show that these can be used to accurately describe the optimal performance that can be achieved  and, in particular, to devise the regimes where the error exponent is an accurate proxy for finite sample-size performances.
 
     \end{enumerate}
 
\subsection{Related works} 

In terms of finite simple-size analysis  within the centralized framework,  \cite{strassen2009asymptotic} presented 
non-asymptotic results for the optimal \textsc{Type II} error probability  under a constant \textsc{Type I} error restriction in the i.i.d case.  Interestingly, the discrepancy between optimal finite-length and asymptotic performance  was characterized, scaling as  $\mathcal{O}(\sqrt{n})$ with the sample size $n$.
In the same communication-free context, \cite{sason2012moderate} borrows ideas from {\em moderate deviation analysis} \cite{dembo2011large} to obtain an interesting upper bound for the Bayesian error probability by bounding the \textsc{Type I}-\textsc{Type II} errors in a way that both decay to zero sub-exponentially with $n$. More recently,  in \cite{espinosa2021finite}, we obtained non-asymptotic upper and lower bounds for the \textsc{Type II} error probability for i.i.d samples draw according to two arbitrary distributions. We showed that the error exponent is a good approximation for the \textsc{Type II} error probability in the finite sample regime. Importantly, the distributed  setting investigated in this work, with a non-trivial rate constraint in one of the modalities,  induces a mathematical problem that is fundamentally different in terms of the requested tools. 

Communication restrictions subject to zero-rate (in bits-per sample) have been investigated in  \cite{watanabe2018neyman}. The error exponent and non-asymptotic bounds have been characterized. Extensions to interactive HT with zero-rate have been reported in~\cite{DBLP:journals/corr/KatzPD16a}.
A preliminary version of this work was presented in \cite{espinosa2019new} with partial results and sketches of some of the arguments. In this paper, we extend the results for a larger family of scenarios, provide complete proofs of the results and more systematic analysis of the practical implications of these results. 

The outline of the paper is as follows. Section \ref{sec_main} introduces the main definitions  and reviews some seminal results for the case of unconstrained communication. Sections \ref{Sec4} and \ref{sec_finite_lenght} present our main theoretical results for the asymptotic and the non-asymptotic regimes, respectively. Numerical analysis and discussions are relegated to Section \ref{sec_numerical}.  Section \ref{sec_final} concludes the paper. Finally,  the proofs are relegated to  Appendix.

\section{Problem Setting and Preliminaries}
\label{sec_main}
Let us consider a finite alphabet product space $\mathbb{Z} = \mathbb{X} \times \mathbb{Y}$, where $\mathcal{P}(\mathbb{Z})$ denotes the family of probabilities on $\mathbb{Z}$. We have a joint random vector $(X,Y)$ with values in $\mathbb{Z}$ and  equipped with a joint probability $P \in \mathcal{P}(\mathbb{Z})$ where $P_X\in \mathcal{P}(\mathbb{X})$ and $P_Y \in \mathcal{P}(\mathbb{Y})$ denote the marginal of  $X$ and $Y$, respectively. $X_1^n=(X_1,...,X_n)$ and $Y_1^n=(Y_1,...,Y_n)$ denote the finite block vector with product  (i.i.d.) distribution $P_{X^n_1 Y^n_1} \triangleq P^n \in \mathcal{P}(\mathbb{X}^n \times \mathbb{Y}^n)$. 
We consider two scenarios for the data generated distribution of $(X^n_1,Y^n_1)$, i.e., 
\begin{equation}\label{eq_sec_main_1}
\begin{array}{ll}
H_0: &   (X_1^n,Y_1^n) \sim  P^n_{XY},\\
H_1: &   (X_1^n,Y_1^n) \sim Q^n_{XY}, 
\end{array} 
\end{equation}
where $Q_{XY} = P_X \cdot P_{Y}$ denote the product probability modeling the case where $X^n_1$ and $Y^n_1$ are independent. In order to make the problem non-trivial, we assume that~\cite{cover2012elements}: 
\begin{align}
\mathcal{D}(P_{XY}\|Q_{XY})&=\sum_{(x,y)\in \mathcal{X}\times \mathcal{Y}} P_{XY}(x,y) \log \frac{ P_{XY}(x,y)}{ Q_{XY}(x,y)}\nonumber\\ &=I(X;Y)>0.
\end{align}
where $\mathcal{D}(\cdot\|\cdot)$ is the divergence between two probabilities and $I(X;Y)$ is the mutual information between $X$ and $Y$~\cite{cover2012elements}. 

Without communication constraints, the fusion center needs to decide about the true underlying hypothesis ($H_0$ or $H_1$) based 
on an observation of the joint vector $(X^n_1,Y^n_1)$. Here we introduce a decentralized version of this problem which is illustrated in Fig. \ref{hypo}. In this distributed context, the decision rule is composed by a pair of encoder and decoder $(f_n,\phi_n)$ of length $n$ and rate $R$ (in bits per sample), where:  
\begin{align}\label{eq_encoder_detector}
&f_n: \mathbb{X}^n \rightarrow \{ 1,\dots,2^{nR} \},\ \textrm{(encoder)}\nonumber\\
&\phi_n: \{1,\dots,2^{nR} \} \times  \mathbb{Y}^n \rightarrow \Theta=\{0,1 \},\ \textrm{(decoder)}.
\end{align}
$f_n(\cdot)$ models a fixed-rate lossy encoder (or quantizer) of $X_1^n$ and $\phi_n(\cdot)$ represents the detector (or classifier) acting on the one-sided compressed data $(f_n(X_1^n),Y_1^n)\in \{ 1,...,2^{nR} \} \times  \mathbb{Y}^n$. 
The encoder represents a remote agent that senses $X_1^n$ and transmit a finite description (using $R$ bits per sample) of $X_1^n$ to a fusion center (see Fig.\ref{hypo}). The fusion center receives the quantization of $X_1^n$ and at the same time senses locally a second modality $Y_1^n$ to guess (using $\phi_n(\cdot)$) about the true distribution of the joint vector $(X^n_1,Y^n_1)$.  For any pair $(f_n,\phi_n)$ of length $n$ and rate $R$, we introduce the corresponding \textsc{Type I} and \textsc{Type II} error probabilities  \cite{kendall1999kendall}, \cite{kullback1951information}:
\begin{align}
\label{eq_sec_main_2}
P_0(f_n,\phi_n)& \triangleq P_{XY}^n\big(\mathcal{A}^c(f_n,\phi_n)\big)\text{ and }\\
\label{eq_sec_main_3}
P_1(f_n,\phi_n)& \triangleq Q_{XY}^n\big(\mathcal{A}(f_n,\phi_n) \big), 
\end{align}
where $A(f_n,\phi_n) \triangleq \{ (x_1^n,y_1^n)\in \mathbb{X}^n \times \mathbb{Y}^n : \phi_n(f_n(x_1^n),y_1^n) = 0\}$.
Traditionally, for any $\epsilon>0$, we are interested in the family of optimal encoder-decoder pairs satisfying: 
\begin{equation} 
\beta_n(\epsilon,R) \triangleq \min_{(f_n,\phi_n)} \big\{P_1(f_n,\phi_n): P_0(f_n,\phi_n)\leq \epsilon \big\},
\label{typeII}
\end{equation}
where the minimum  is over all encoding-decoder pairs in (\ref{eq_encoder_detector}).

We study the performance of the optimal scheme (\ref{typeII}) by focusing on the case where a sequence of restrictions $(\epsilon_n)_{n\geq 1}$ is required to tend to zero as the sample size grows. The objective is to explore how this restriction is expressed in the terms of $(\beta_n(\epsilon_n,R))_n$ with $n$ in conjunction with other properties of the problem (e.g., the distribution $P_{XY}$ and the rate $R$). In this work, we are primarily interested in deriving expressions to bound $\beta_n(\epsilon_n,R)$ in the large sample regime (non-asymptotic). To this end, it is essential to first characterize the asymptotic nature of the sequence $(\beta_n(\epsilon_n,R))_n$. However, before presenting the main contributions of this paper, we review  some essential asymptotic results for the classical communication-free (centralized) scenario.  


\subsection{Review of centralized HT results}
For completeness, it is worth revisiting the centralized case where $f_n: \mathbb{X}^n \rightarrow \mathbb{X}^n$ is the identity  
and the solution of (\ref{typeII}) is then denoted by $\beta_n(\epsilon_n)$. Furthermore, when $\epsilon_n=\epsilon>0$ for all $n$,  
this is a classical HT setting where the celebrated Stein's Lemma implies the following result~\cite{chernoff1952measure,cover2012elements}:
\begin{lemma}[\textit{Stein's Lemma}]
	\label{lemma_stein}
	For any $\epsilon \in (0,1)$, 
	$$
	\lim_{n \rightarrow \infty} -\frac{1}{n } \log \beta_n(\epsilon) =  \mathcal{D}(P\|Q).
	$$
\end{lemma}
This result establishes the asymptotic decayment of the \textsc{Type II} error subject to a fixed $\epsilon>0$ implying that $\beta_n(\epsilon) \sim e^{-n \mathcal{D}(P\|Q)}$ as $n$ tends to infinity (large sampling regime). Interestingly, in \cite{espinosa2021finite}, we provided upper and lower bounds in the finite length regime showing that in practice the number of samples required to approximate the \textsc{Type II} error probability  to $(e^{-n \mathcal{D}(P\|Q)})$ is not large. This observation supports the claim that the exponential approximation is a useful proxy for \textsc{Type II} error probability. 


As a matter of fact, for the sub-exponential regime of $(\epsilon_n)_n$, the following result is known. 
\begin{lemma}
\label{lemma_stein_restr}(\cite[Sect. IX]{nakagawa1993converse})
if $(1/\epsilon_n)$ is $o(e^{rn})$ for any $r >0$ then $\lim_{n\rightarrow \infty} - \frac{1}{n} \log \beta_n(\epsilon_n)=\mathcal{D}(P\|Q)$.
\end{lemma}

Therefore, the error exponent obtained with a fixed $\epsilon>0$ in Lemma \ref{lemma_stein}
is preserved for a family of stringer decision problems in (\ref{typeII}) as long as $({\epsilon_n})_n$ tends to zero 
at a sub-exponential rate.

\subsection{Review of distributed HT results}
Returning to the main decentralized task with communication constraints in Fig.\ref{hypo},   \cite{ahlswede1986hypothesis} determined the following result\footnote{This result can be interpreted as the counterpart of the Stein's Lemma in the  decentralized setting of Fig.\ref{hypo}.}: 
\begin{lemma} \label{th_ahlswede&csiszarlimit}
 \cite[Theorem 3]{ahlswede1986hypothesis} For any $\epsilon >0$, it follows that\footnote{This result provides an interesting connection with the problem 
 noisy lossy source coding  with  log-loss fidelity  \cite{shkel}.  The performance limits in the right hand side (RHS) of (\ref{csiszarlimit}) coincides precisely with the distortion-rate function of the information bottleneck problem~\cite{tishby2000information}.}
\begin{equation}
\xi(R) \triangleq  \lim_{n\rightarrow \infty}  - \frac{1}{n}\log  \beta_n(\epsilon,R) =  \max_{\begin{subarray}{c} U : U \mkv X \mkv Y \\ 
I(U;X) \leq R \ |\mathbb{U}| \leq |\mathbb{X}|+1
\end{subarray}}I(U; Y),
\label{csiszarlimit}
\end{equation}
where $ U \mkv X \mkv Y$ denotes the fact that $(U,X,Y)$ forms a Markov chain (i.e., $(U,Y)$ are independent condition to $X$). 
\end{lemma}
The result presented in (\ref{csiszarlimit}) is a trade-off between representation and regularisation, in the sense that we seek to learn the best possible representation of $X$ for predicting $Y$. As for the more challenging scenario where $(\epsilon_n)_n$ tends to zero with $n$, in \cite{han1987hypothesis} the author provided a lower bound for the error exponent of the \textsc{Type II} error probability  
in the case of exponentially decreasing \textsc{Type I} error restrictions:
\begin{lemma}\cite[Han and Kobayashi]{han1989exponential}
\label{Hanmain}
Let us assume that $\epsilon_n \leq e^{-rn}$ for some $r>0$, then: $\liminf\limits_{n\rightarrow \infty} -\frac{1}{n}\log  \beta_n(\epsilon_n,R) \geq$
\begin{align}
&\max_{w \in \rho(R,r)} \min_{\begin{subarray}{c} \tilde{P}_{UXY} \\ 
\mathcal{D}(\tilde{P}_{UXY}\|P_{UXY})\leq r \\
\tilde{P}_{U|X}=P_{U|X}=w \\
U \mkv X\mkv Y
\end{subarray}}[\mathcal{D}(\tilde{P}_{X}\|P_{X})+I(U;Y)],
\label{Hanmain2}\\
&\rho(R,r)\triangleq \{ w \in \mathcal{P}(\mathbb{U}|\mathbb{X})| \max_{\begin{subarray}{c}\tilde{P}_{X}: \mathcal{D}(\tilde{P}_{X}\| Q_X)\leq r \\ \tilde{P}_{U|X}=w\\ P_{UX}=w\cdot\tilde{P}_{X} \end{subarray}}\limits I(U;X) \leq R \},\nonumber 
\end{align}
$\mathcal{P}(\mathbb{U}|\mathbb{X})$ denotes all test (quantizer) channels from $\mathbb{X}$ to $\mathbb{U}$.
\end{lemma}

\section{Asymptotic Result}
\label{Sec4}
Our first result complement the regime on $(\epsilon_n)_n$ presented in Lemma \ref{Hanmain} to obtain an asymptotic characterization of $(\beta_n(\epsilon_n,R))_n$. 
In particular,  we explore the 
important sub-exponential regime for the restriction sequence  $(\epsilon_n)_n$ of \textsc{Type I} error probability. The proof is relegated to Appendix~ \ref{proofresult}.
\begin{theorem}
\label{overheadR}
Let us assume that $(1/{\epsilon_n})_n=o(e^{rn})$ for any $r >0$. Then, 
\begin{equation}
\lim_{n\rightarrow \infty}\limits - \frac{1}{n}\log (\beta_n(\epsilon_n,R))= \xi(R),
\label{achievableR}
\end{equation}
where $\xi(R)$ is defined in (\ref{csiszarlimit}).
\end{theorem}
This result establishes an extensive regime on the velocity at which $({\epsilon_n})_n$ tends to zero for which the error exponent of the problem is invariant and matches the expression obtained for the less restrictive and classical setting ($\epsilon_n=\epsilon>0$) 
presented in Lemma \ref{th_ahlswede&csiszarlimit}.
This result is interesting because, as it was pointed out in \cite{han1989exponential}, there was no guarantee that the asymptotic limit in (\ref{achievableR}) remains the same as the result in Lemma \ref{th_ahlswede&csiszarlimit} when moving to stringer regimes on the velocity at which $({\epsilon_n})_n$ vanishes with $n$.  
Besides this, the above result can be considered as being the counterpart of what is observed in the centralized setup when contrasting   Lemmas \ref{lemma_stein_restr} and \ref{lemma_stein}.

The proof of Theorem \ref{overheadR}, relegated to  Appendix \ref{proof_overheadR}, is divided into two parts. 
The direct part (i.e., constructive argument) is based on constructing an encoder-decision pair that guarantees that the error exponent of the optimal \textsc{Type II} is greater than $\xi(R)$. The second part of the argument (i.e., the infeasibility part) 
proves that no pair of encoder-decoder rule satisfying the restriction of the  \textsc{Type I} error has an error exponent greater than $\xi(R)$.
%
The proof argument used in both the achievable and infeasibility parts (see Appendix \ref{proof_overheadR}) 
is based on a refined use of concentration inequalities \cite{boucheron2013concentration}. 
In particular, 
the achievable part is divided into two steps. The first step consists of reducing the problem to an i.i.d. structure over a block of $X_1^n$ induced by the encoder, which will concentrate (in probability) to an  error exponent that is different from $\xi(R)$ in (\ref{achievableR}). Importantly, the discrepancy between the concentration limit obtained from our approach (i.e., finite-block strategy) and $\xi(R)$ can be resolved analytically by connecting our problem with a noisy rate-distortion problem, where the discrepancy between its fundamental limit and a finite length version of this object is well understood \cite{zhang1997redundancy}. The second step consists of optimizing our approach by giving concrete conditions to make the discrepancy between $\xi(R)$ and $- \frac{1}{n}\log (\beta_n(\epsilon_n,R))$ vanishes with $n$.

\section{Finite-length Result} 
\label{sec_finite_lenght}
Our main result is concerned with  the practically relevant task of offering a non-asymptotic characterization of the sequence $(\beta_n(\epsilon_n,R))_n$ for different scenarios of $(\epsilon_n)_n$, given the model $P_{XY}$ and the rate constraint $R>0$. To address this question, our methodology  uses the asymptotic limit of $(\beta_n(\epsilon_n,R))_n$, stated in Theorem \ref{overheadR}, and from this, analizes the discrepancy between $- \frac{1}{n}\log \beta_n(\epsilon_n,R)$ and $\xi(R)$ as a function of $n$. In concrete, our main result (stated below) derives upper and lower bounds for  $- \frac{1}{n}\log \beta_n(\epsilon_n,R)$ in different sub-exponential scenarios for 
the  \textsc{Type I} restriction sequence $(\epsilon_n)_n$. As a  corollary,   we determine the velocity at which  $- \frac{1}{n}\log \beta_n(\epsilon_n,R)$ achieves its limit in (\ref{achievableR}). The proof is relegated to Appendix \ref{proofresult2}.

\begin{theorem}\label{theorem2}
Assume that  
$R<H(X)$. Then, 
\begin{itemize}
\item[i)] If $(\epsilon_n)_n =(1/\log(n))_n$ (logarithmic), it follows: 
\begin{align}
  -\frac{1}{n}&\log(\beta_n(\epsilon_n,R)) -\xi(R) \geq \nonumber\\
  & \left( \frac{\partial D(R)}{6\partial R}  -\frac{\sqrt{2\ln(\log(n))}C(P_{XY})}{\log(n)}- o\left(1\right )\right )  \frac{\log n}{n^{1/3}} \label{converse1}
  \end{align}
  \begin{align}
    -\frac{1}{n}&\log(\beta_n(\epsilon_n,R)) -\xi(R)   \nonumber \\
   &\leq   \left (16C(P_{XY})+\frac{\log(\log(n))\sqrt{\log(n)}}{n} \right )  \frac{1}{\sqrt{\log(n)}};
\label{overheadtheoremv3}
\end{align} 

\item[ii)] If $(\epsilon_n)_n =(1/n^{p})_n$ (polynomial) with $2>p>0$, then
\begin{align}
 -\frac{1}{n}&\log(\beta_n(\epsilon_n,R)) -\xi(R)\geq \nonumber \\
 & \left ( \frac{1}{6} \frac{\partial D(R)}{\partial R} -\frac{\sqrt{2p\ln(n)}}{\log n}  C(P_{XY}) - o\left(1\right )\right )  \frac{\log n}{n^{1/3}}  \label{converse2} \\
-\frac{1}{n}&\log(\beta_n(\epsilon_n,R)) -\xi(R) \nonumber \\
&\leq \left (16C(P_{XY})+\frac{p\log(n)}{n^{1-p/2}} \right ) \frac{1}{n^{p/2}};
\label{overheadtheorem2}
\end{align} 

\item[iii)] If $(\epsilon_n)_n =(1/n^{p})_n$ (polynomial) with $p\geq 2$, then
\begin{align}
 -\frac{1}{n}&\log(\beta_n(\epsilon_n,R))-\xi(R)\geq \nonumber \\ 
  &  \left ( \frac{1}{6} \frac{\partial D(R)}{\partial R} -\frac{\sqrt{2p\ln(n)}}{\log n}  C(P_{XY}) - o\left(1\right )\right )  \frac{\log n}{n^{1/3}} \label{converse3}  \\
 -\frac{1}{n}&\log(\beta_n(\epsilon_n,R)) -\xi(R) \nonumber \\ 
 & \leq  \left (8\sqrt{2}C(P_{XY})\frac{\sqrt{n^{2-p}+1}}{\log(n)}+2\right )  \frac{\log(n)}{n};
\label{overheadtheoremv4}
\end{align} 

\item[iv)] If $(\epsilon_n)_n =(1/e^{n^{p}})_n$ (superpolynomial) with $p\in (0,1 )$, 
\begin{align}
-\frac{1}{n}&\log(\beta_n(\epsilon_n,R))-\xi(R) \geq \nonumber \\
 & \left ( \frac{(1-p)}{6} \frac{\partial D(R)}{\partial R} -\frac{\sqrt{2} C(P_{XY})}{\log(n)}- o\left( 1\right )\right )   \frac{\log n}{n^{(1-p)/3}}\label{converse4} \\
  -\frac{1}{n}&\log(\beta_n(\epsilon_n,R))-\xi(R) \nonumber \\
  &\leq \left (8\sqrt{2}C(P_{XY})\frac{\sqrt{e^{-n^{p}}n^2+1}}{\log(n)}+2\right ) \frac{\log(n)}{n}.
\label{overheadtheoremv2}
\end{align} 

$D(R)$ is the noisy distortion-rate function \cite{cover2012elements} and \\
$C(P_{XY})\triangleq \sup_{(x,y) \in \mathbb{X} \times \mathbb{Y}}\limits \Big|  \log\left( \frac{P_{XY}(\{ (x,y) \})}{Q_{XY}(\{ (x,y) \})}\right)\Big|  < \infty
$.
\end{itemize}
\end{theorem}

\subsection{Discussion of Theorem \ref{theorem2}}
\label{subsec_discussion_th2}

 \textbf{(i)} The results establish non-asymptotic  bounds for the \textsc{Type II} error when we impose concrete scenarios for the monotonic behavior of $(\epsilon_n)_n$. We explore three main regimes for $(\epsilon_n)_n$: logarithmic, polynomial, super-polynomial. Each of these cases has its corresponding lower and upper bounds,  which depends specifically on the considered  $(\epsilon_n)_n$. 

 \textbf{(ii)} The proof of Theorem \ref{theorem2} involves an optimization problem of the upper and lower bounds presented in the proof of Theorem \ref{overheadR}, for which the arguments used to prove Theorem \ref{overheadR} were instrumental for this analysis. Specifically, we refine the analysis introduced in (\ref{overhead_final}), (\ref{newdelta}) and (\ref{eq_final_equation}) by finding optimal values for $l$ and $s_n$ for a given $\epsilon_n$. These choices of values for $l$ and $s_n$ give us non-asymptotic lower and upper bounds for $-\frac{1}{n}\log(\beta_n(\epsilon_n,R))$, for each scenario.

\textbf{(iii)} Regarding the upper bound of  $-\frac{1}{n}\log(\beta_n(\epsilon_n,R))$ ((\ref{overheadtheoremv3}), (\ref{overheadtheorem2}), (\ref{overheadtheoremv4}) and (\ref{overheadtheoremv2})), obtained from the impossibility argument (converse part), as $({\epsilon_n})_n$ goes  to zero faster (from case to case),  the velocity at which 
the bound tends to zero increases; from the slower rate $\mathcal{O}\left ({1}/{\sqrt{\log(n)}}\right)$ to the faster that is $\mathcal{O}\left ({\log(n)}/{n} \right )$. Therefore, by imposing a more restrictive $({\epsilon_n})_n$ there is an effect in the discrepancy between the fundamental limit $\xi(R)$ and the optimal \textsc{Type II} error $-\frac{1}{n}\log \beta_n(\epsilon_n,R)$ obtained from this upper bound analysis.

\textbf{(iv)} Regarding the lower bound of $-\frac{1}{n}\log \beta_n(\epsilon_n,R)$ ((\ref{converse1}), (\ref{converse2}), (\ref{converse3}) and (\ref{converse4})), obtained from the direct argument (achievability part), as $({\epsilon_n})_n$ goes faster to zero (from case to case), the derived bound -for the super-polynomial case- decreases in the velocity at which the discrepancy in error exponent (i.e., $-\frac{1}{n}\log(\beta_n(\epsilon_n,R)) -\xi(R)$) tends to zero.  For the other two cases (logarithmic and polynomial), the velocity is not affected, but the constants change to slower magnitudes. These trends are consistent with the observation that by relaxing the velocity of $({\epsilon_n})_n$ the decision problem is less restrictive and then, the result favors the possibility of obtaining a better \textsc{Type II} error (smaller) than the one predicted by the asymptotic limit, which is $e^{-n\xi(R)}$. 

\textbf{(v)} Finally, it is worth noting that if we consider the relaxed restriction $\epsilon_n=\epsilon \in (0,1)$ in Lemma \ref{th_ahlswede&csiszarlimit}, the achievability part of our argument still works and for $\xi(R)- \left ( - \frac{1}{n}\log \beta_n(\epsilon,R) \right )$ it offers  an upper bound that converges to zero as $\mathcal{O}\left( \frac{ \log(n)}{n^{1/3}}\right )$. 

This last velocity of convergence is slower than the same result known for the unconstrained (centralized)  problem presented in \cite{strassen2009asymptotic}.
In fact, when $X^n_1$ is fully observed at the detector (see  Lemma \ref{lemma_stein}), in \cite{strassen2009asymptotic} the author showed that the discrepancy $\left \rvert \mathcal{D}(P\|Q)- \left( -\frac{1}{n}\log \beta_n(\epsilon) \right )\right \rvert$ 
tends to zero as $\mathcal{O}\left({1}/{\sqrt{n}}\right)$.\footnote{For completeness, this is presented in Lemma \ref{th_strassen2009asymptotic} in Appendix \ref{app_finite-length_unconstrained}}. 
We conjecture that our slower rate can be attributed to the non-trivial role of the communication constraint in our problem, 
which  breaks the i.i.d. structure of $X^n_1$  in a way that it is not possible to use the tools adopted to 
derive the unconstrained result in Lemma \ref{th_strassen2009asymptotic}. 
It is a topic of further research to uncover if the upper bound $\mathcal{O}\left( \frac{ \log(n)}{n^{1/3}}\right )$ for the discrepancy $\xi(R)- \left ( - \frac{1}{n}\log \beta_n(\epsilon,R) \right )$  can be improved, or if it is possible  to show (by a converse argument) that this rate is indeed optimal provided that  $\epsilon_n=\epsilon>0$.

\subsection{Interpretation of Theorem \ref{theorem2}}
\label{sub_sec_inter_nun}
In general, Theorem \ref{theorem2} can be presented as two bounds: 
\begin{align}\label{eq_sec_inter_nun_1a}
	\xi_o - f(n) \leq \frac{1}{n} &\log \beta_n  \leq \xi_o + g(n), 
\end{align}
where $\beta_n$ is the optimal \textsc{Type II} error consistent with \textsc{Type I} error restriction ($\epsilon_n$ in the statement of Theorem \ref{theorem2}), $\xi_o$ is the performance limit (in Theorem \ref{overheadR}),  $f(n)$ is a positive sequences  that goes to zero with $n$ ($o(1)$) representing the penalization (in error exponent)  for the use of finite simple-size, and $g(n)$ is a positive sequence that goes to zero representing a discrepancy with the limit but that can be seen as a gain in error exponent.  
%
Then, we have a feasibility range for $\beta_n$ given by the interval: 
$$
\big[\exp[-n( \xi_o+g(n))], \exp[ -n (\xi_o-f(n))]\big].
$$
This interval contains the nominal value $ e^{ -n \xi_o} $, which is consistent with the error exponent limit in Theorem \ref{overheadR} but extrapolated to a finite length regime.  If we consider $ \exp (-n \xi_o) $ as our reference,  we can study two feasible regions: the pessimistic interval  $(\exp( -n \xi_o), \exp( -n (\xi_o-f(n)))]$ where the error probability is greater than the nominal value $e^{ -n \xi_o} $, and the optimistic interval  $[\exp(-n (\xi_o+ g(n))), \exp( -n \xi_o)]$ where the appositive occurs. The length of the interval of the two regions  is an  indicator of the precision of our result  (the worse case discrepancy with respect to 
$e^{ -n \xi_o} $). For the pessimistic region,  the length of that interval is $e^{-n\xi_o}(e^{nf(n)}-1)$.  
From the fact that  $f(n)$ is $o(1)$ (see the statement of Theorem \ref{theorem2}), the length of this interval tends to zero strictly faster  than $\mathcal{O}(e^{ -n (\xi_o-\epsilon)})$ for any $\epsilon>0$ and, consequently, the precision has an exponential rate of convergence  that is asymptotically given by the nominal exponent $\xi_o>0$. On the optimistic region, the length of this interval is $e^{-n\xi_o}(1-e^{-n g(n)})$, which is  $\mathcal{O}(e^{-n \xi_o})$. Overall, the length of the pessimistic interval dominates the analysis and, consequently,  the precision of the result (i.e., the worse case discrepancy with respect to the nominal $ e^{ -n \xi_o} $) tends to zero as  $\mathcal{O}(e^{ -n (\xi_o- f(n))})$. This order  is equivalent  to the worse-case \textsc{Type II} error probability  ($e^{ -n( \xi_o- f(n))}$) predicted from Theorem \ref{theorem2}. 

In conclusion, the overall quality of the result is governed  by $\xi_o$ and affected in a smaller degree by how fast $f(n)$ goes to zero.  Note that  $g(n)$ plays no role from this perspective.  We discussed on the previous section that $f(n)$ goes faster to zero when we relax the problem (i.e., passing from a scenario for $({\epsilon_n})_n$ to a scenario where this sequence tends to zero at a smaller  velocity). Then, the precision of Theorem \ref{theorem2}  improves when simplifying the problem from one restriction  $({\epsilon_n})_n$ to a relaxed restriction $(\tilde{\epsilon}_n)_n$ for the \textsc{Type I} error. This reinforces one of the points mentioned in Section \ref{subsec_discussion_th2}, where we discussed that the velocity at which $({\epsilon_n})_n$ goes to zero does not affect the limit $\xi_o$ (Theorem \ref{overheadR}) but it does affect our finite length result 
through $f(n)$. 

\section{Application Examples}
\label{sec_numerical}
In this section, we present some empirical evidences illustrating the possible implication  of Theorem \ref{theorem2} to effectively bound $\beta_n(\epsilon_n,R)$ with finite-sample size $n$. 
Theorem \ref{theorem2} offers an interval of feasibility for $\beta_n(\epsilon_n,R)$ 
expressed by
\begin{align}
 &\textrm{UB}(\epsilon_n,R)=\exp\left[-n\left(\xi(R)+\frac{\partial D(R)}{\partial R}\frac{\log (l)}{2l}\right.\right. \nonumber \\
&\left.\left. -\sqrt{\frac{2l\ln(1/\epsilon_n)}{n}} C(P_{XY})\right)\right],\label{Achievable_bound}\\
& \textrm{LB}(\epsilon_n,R)= \exp\left[-n\left(\xi(R)+4C(P_{XY})\cdot \right.\right. \nonumber \\
&\sqrt{2\ln\left( \frac{1}{1-\epsilon_n-h_n(s)}\right)}+\left.\left.\frac{\log(1/h_n(s))}{n}\right)\right],
\label{Converse_bound}
\end{align}
where $\beta_n(\epsilon_n,R) \in [  \textrm{LB}(\epsilon_n,R),  \textrm{UB}(\epsilon_n,R)]$.\footnote{$l$ and $h_n(s)$ are obtained according to the proof of Theorems \ref{overheadR} and \ref{theorem2} (see Appendix~\ref{proofresult} and \ref{proofresult2} for details).}  
The length of $[\textrm{LB}(\epsilon_n,R),\textrm{UB}(\epsilon_n,R)]$ indicates the precision of our approximation and the interval itself can be used to evaluate how 
representative  is $e^{-n \xi(R)}$ of $\beta_n(\epsilon_n,R)$ for a finite $n$.

We first evaluate the length of $[\textrm{LB}(\epsilon_n,R),\textrm{UB}(\epsilon_n,R)]$ by 
considering  
four cases $(\epsilon_n)_n \in \{0.01,1/\log(n),n^{-0.01},n^{-0.1}$ $ \}$ associated to a constant, a logarithmic and a polynomial \textsc{Type I} error restriction, respectively. We use a discretized version of a Gaussian pdf $P_{XY}$ 
of $|\mathbb{X}|\times |\mathbb{Y}|$ where the mutual information between the two variables ($X$ and $Y$) is 7 and 1.5 nats, respectively.  
To compute the expressions in (\ref{Achievable_bound}) and (\ref{Converse_bound}), we need to evaluate $\xi(R)$.  Obtaining $\xi(R)$ involves an optimization problem with respect to the encoder $f_n$ and the rate $R$ \cite{tishby2000information}. To this end, 
we use the algorithm in \cite{vera2018compression} which is a generalization of {\em Blahut-Arimoto algorithm} \cite{arimoto1972algorithm}.\footnote{Importantly, under some mild conditions given in \cite{vera2018compression}, this optimization (algorithm) converges to $\xi(R)$.} 

\begin{table*}[ht]
\centering
\begin{tabular}{|c|c|c|c|c|c|c|c|c|}
\hline
  & \multicolumn{8}{c|}{Number of observations $n$ } \\
\hline
 $\epsilon_n$     & 50 & 150 & 250 & 350 & 450 & 550 & 650 & 750  \\
\hline
 $1/\log(n)$  &       1.2138e-12&  2.7758e-37 &  3.3636e-62 & 3.0109e-87& 2.2227e-112& 1.4286e-137& 8.2535e-163 & 4.3764e-188 \\
$1/n^{0.01}$  &     6.4432e-10 & 1.4348e-30& 4.2359e-52&  4.5381e-74&  2.5069e-96& 8.5307e-119& 1.9894e-141&3.4117e-164\\
$1/n^{0.1}$   &          0.0045 & 1.8737e-14 &  2.5598e-28&  2.0497e-43 & 2.5977e-59 & 8.5949e-76 & 1.0006e-92 &4.9903e-110 \\
\hline
\end{tabular}
\vspace{2mm}
\caption{Magnitude of $  \textrm{UB}(\epsilon_n)- \textrm{LB}(\epsilon_n)$ function of $\epsilon_n$ and $n$ for the case when $I(X;Y)=1$.}
\label{Tablap}
\end{table*}

\begin{figure}
    \centering
    \begin{minipage}{0.5\textwidth}
        \centering
        \includegraphics[width=1.0\linewidth]{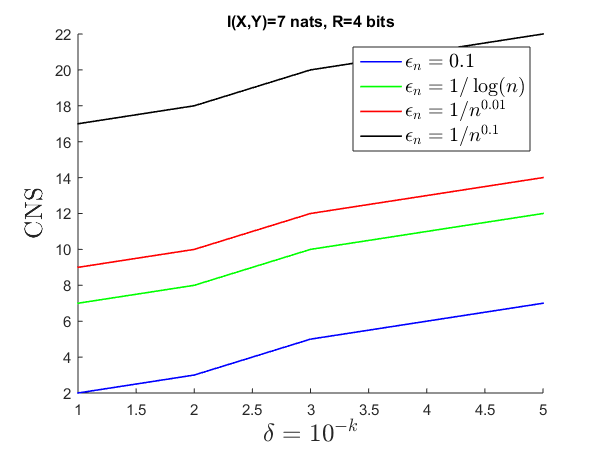} 
        \caption{Critical Number of Samples (CNS) predicted by Theorem~\ref{theorem2}  
        across different values of $\delta=10^{-k}$. 
        The values used are $\xi(R)=3$, $I(X;Y)=7$, $R=4$ and  $C_X(P,Q)=2.47$. 
        }
        \label{High}
    \end{minipage}\hfill
    \begin{minipage}{0.5\textwidth}
        \centering
       \includegraphics[width=1.0\linewidth]{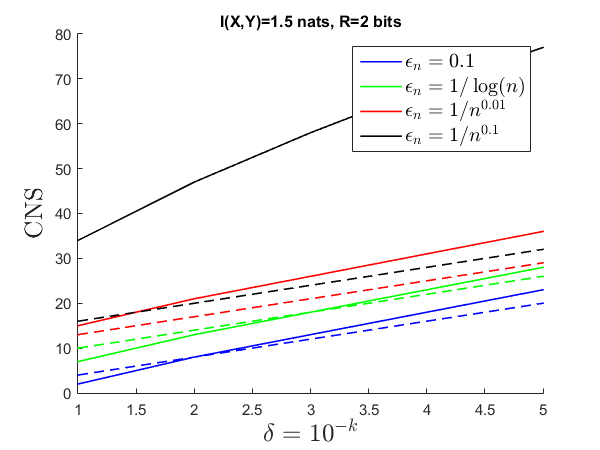} 
        \caption{
        CNS predicted by Th. \ref{theorem2} 
        across different values of $\delta=10^{-k}$. 
        Low rate case with $\xi(R)=0.7$, $I(X;Y)=1.5$, $R=2$ and  $C_X(P,Q)=1.92$. 
        The dashed lines show an estimation of  the exact  CNS obtained from $\beta_n(\epsilon_n,R)$. 
        }
           \label{Low}
    \end{minipage}
\end{figure}

Table \ref{Tablap} shows the lengths of $[\textrm{LB}(\epsilon_n,R),\textrm{UB}(\epsilon_n,R)]$.
We verify that $\textrm{UB}(\epsilon_n,R)-\textrm{LB}(\epsilon_n,R)$ tends to zero exponentially fast with the sample size as observed in Section \ref{sub_sec_inter_nun}. From this exponential decay,  the nominal value predicted by Theorem \ref{overheadR}, i.e.,  $\exp(-n \xi(R))$,  is a very precise approximation of $\beta_n(\epsilon_n,R)$ provided that $n$ is sufficiently large. This support the idea that $e^{-n \xi(R)}$ is an excellent proxy of $\beta_n(\epsilon_n,R)$ if a critical number of samples is achieved. 
Table \ref{Tablap} also shows that the precision of the result measured by $(\textrm{UB}(\epsilon_n,R)-\textrm{LB}(\epsilon_n,R))$ is affected by the velocity at which the \textsc{Type I} error sequence tends to zero, which is consistent 
with our previous analysis in Section \ref{sub_sec_inter_nun}.
In particular, we observe that for a faster convergence rate of $({\epsilon_n})_n$, i.e, a stringer distributed decision problem, the length of $[\textrm{LB}(\epsilon_n,R),\textrm{UB}(\epsilon_n,R)]$ is bigger, which means that our bounds are expected to be less informative on $\beta_n(\epsilon_n,R)$ when compared with a relaxed  scenario.

The results  presented in Table \ref{Tablap}  support the claim that $\exp(-n \xi(R))$ can be adopted as practical proxy to  
$\beta_n(\epsilon_n,R)$. To formalize this, we address the following question:  
for a given arbitrary small $\delta >0$ of the form $10^{-k}$ with $k \in \left\{1,..,5 \right\}$ and a joint model $P_{XY}$, we seek to find the lowest $n$ such that $\beta_n(\epsilon_n,R) \in (e^{-n \xi(R)}- \delta, e^{-n \xi(R)} + \delta)$. The exponential decay of the length of $[\textrm{LB}(\epsilon_n,R),\textrm{UB}(\epsilon_n,R)]$,  observed in Table  \ref{Tablap}, suggests that this condition happens eventually with $n$ very quickly. Importantly, we can derive an upper bound for this Critical Number of Samples (CNS) from the closed-form expressions we have for $\textrm{LB}(\epsilon_n,R)$ and $\textrm{UB}(\epsilon_n,R)$.\footnote{The predicted CNS is the first $n\geq 1$ such that $\max\{ \textrm{UB}(\epsilon_n,R)- e^{-n \xi(R)},e^{-n \xi(R)}-\textrm{LB}(\epsilon_n,R) \}  \leq \delta$, which is finite for any $\delta>0$ and can be computed from our result.} Figs. \ref{High} and \ref{Low} present the predicted CNS vs. $\delta= 10^{-k}$ for different scenarios of $P_{XY}$ (in terms of the magnitude of $I(X;Y)$) and $(\epsilon_n)_n$. We consider two scenarios for $P_{XY}$ 
($I(X;Y)=7$, $R=4$ and $I(X;Y)=1.5$ with $R=2$) and we explore $(\epsilon_n)_n  \in \{n^{-0.01}, n^{-0.1},1/\log(n), 0.1\}$. Figs. \ref{High} and \ref{Low} show that even for a very small precision $\delta=10^{-5}$, the point at which $\beta_n(\epsilon_n,R)$ is well approximated by $e^{-n\xi(R)}$ happens with less than $22$ samples for the high-rate restriction case and in less than 80 samples for the low rate case for the majority of $(\epsilon_n)_n$.\footnote{The observed variations can be attributed to the value of $\frac{\partial D(R)}{\partial R}$, which tends to zero as long as $R>H(X)$.}
The dependency of these values (predicted CNS from Theorem \ref{theorem2}) on the magnitude of $I(X;Y)$ and 
$(\epsilon_n)_n$ is clearly expressed, which is consistent with our previous analyses in Section \ref{sub_sec_inter_nun}.

Finally, to evaluate the tightness of our theoretical bounds for the CNS, we simulate data from the true model $P_{XY}$ 
(i.i.d. samples) to have a practical lower bound for $\beta_n(\epsilon_n)$. In particular, given $P_{XY}$, $R$ and $(\epsilon_n)_n$, we obtained empirical estimations of the two error probabilities from which we estimate $\beta_n(\epsilon_n,R)$. $2.5\cdot 10^6$ realizations of $P_{XY}$  were used to obtain good estimations of these probabilities.\footnote{To achieve this, we use an scalar quantization based on the {\em Lloyd-max algorithm} \cite{gallager2008principles} to obtain an induced quantized distribution $P_{f(X_1^n)Y}$.} Using the estimated values of $\beta_n(\epsilon_n,R)$, we obtained for each $\delta >0$ the corresponding CNS where the condition $\beta_n(\epsilon_n,R) \in (e^{-n \xi(R)}- \delta, e^{-n \xi(R)} + \delta)$ is meet directly (the empirical estimations of the CNS). Fig. \ref{Low} contrasts our predictions (theoretical upper bounds) with the empirical estimations  (the dashed lines) of the CNS.
Consistent with the nature of our result, the predicted CNS values are more conservative than the CNS estimated from simulations. 
This discrepancy is not significant overall, in particular for the regime when $\epsilon_n$ exhibits a relatively small velocity of convergence to zero.
Overall, we can conclude that the derived bounds are meaningful   and can be adopted in cases where it is prohibitive  to estimate $\beta_n(\epsilon_n,R)$ from data. Indeed, we face this issue in this analysis, as it was not possible to estimate  $\beta_n(\epsilon_n,R)$ for the higher rate cases.\footnote{$\beta_n(\epsilon_n,R)$ is of  order $O(e^{- n \xi(R)})$ so when $R$ is relatively high, the value of $\xi(R)$ tends to $I(X;Y)$ for which  $e^{n I(X;Y)}$ simulations are needed. This number becomes prohibitive, even for $n$ of order of $30$ when $I(X;Y)=1.5$.}
\section{Summary and Concluding Remarks}
\label{sec_final}
This paper explores the problem of testing against independence with one-sided  communication constraints. More specifically, the scenario of two memoryless sources is considered where one of the modalities is transmitted to the decision-maker (fusion center) over a rate-limited channel.  In this context, we explored a general family of optimal tests (in the sense of Neyman-Pearson) where restrictions on the \textsc{Type I} error are imposed.  We are interested in the velocity at which the \textsc{Type II} error vanishes with the sample size. From a  theoretical perspective,  we obtained  the performance limits for a rich family of problems with a decreasing  sequence of  \textsc{Type I} error probabilities  (Theorem \ref{overheadR}). This result stipulates that the error exponent of the \textsc{Type II} error probability tends to an error exponent (fundamental limit) in the form of the classical {\em Stein's Lemma}. This error exponent is expressed in a closed-form, which is a function of the operational rate (in bits per sample) imposed on one of the information sources. Interestingly, this result  implies  that for a large family of \textsc{Type I} error restrictions (vanishing to zero with the sample size), the error exponent coincides with the result obtained in the (classical) scenario where the \textsc{Type I} error restriction is  constant with $n$ (Lemma \ref{th_ahlswede&csiszarlimit}).

Concerning the finite-sample size analysis, our main result (Theorem \ref{theorem2}) provides bounds for the \textsc{Type II} error probability. Using results from rate-distortion theory and concentration inequalities, we obtained upper and lower bounds for this error as a function of $n$ (the number samples),  the sequence $({\epsilon_n})_n$ that models the restriction for the \textsc{Type I} error probability and the underlying distributions. We observed that the bounds offer an interval of feasibility for the optimal \textsc{Type II} error probability, which presents an accurate description.  A closed-form expression for the worse-case \textsc{Type II} error probability was derived  where a discrepancy in the error exponent (with respect to the asymptotic exponent) was identified. This discrepancy (overhead) can be attributed to using a finite number of samples in the decision.  Furthermore, this penalization vanishes at a velocity that is a function of $({\epsilon_n})_n$, and consequently,  we observed the effect of the \textsc{Type I} error restriction in this non-asymptotic analysis.

Finally, we observed that the \textsc{Type II} error probability is arbitrary close (with $n$) to the nominal value predicted by the asymptotic result  $e^{-n\xi(R)}$, where  $\xi(R)$ is the limit in Theorem \ref{overheadR}. Furthermore, the precision in Theorem \ref{theorem2}, measured by the length of the feasible interval, tends to zero exponentially fast. Numerical analysis in some concrete scenarios confirms the predicted quality of the non-asymptotic results in Theorem \ref{theorem2}.

\appendix 
\label{proofresult}

\subsection{Proof of Theorem \ref{overheadR}:}
\label{proof_overheadR}
The proof  is divided in two parts: a lower and an upper bound result.
We begin with the following bound that extend the result presented in  \cite[Theorem 3]{ahlswede1986hypothesis}.

\begin{theorem}
\label{lemma_achievability}
Let us assume that $\epsilon_n>0$ for all $n$ and $(1/{\epsilon_n})_n=o(e^{rn})$ for any $r >0$, then
\begin{equation}
\liminf_{n\rightarrow \infty} - \frac{1}{n}\log(\beta_n(\epsilon_n,R))\geq \xi(R).
\end{equation}
\end{theorem}
\label{sub_proof_1}
\begin{IEEEproof}
For an arbitrary encoder $f_n: \mathbb{X}^n \mapsto \{1,\dots,2^{nR} \} $ of rate $R>0$, let us consider the corresponding optimal decision regions -according to Neyman-Pearson's Lemma- on the one-sided quantized space $\{1,\dots,2^{nR} \} 
\times \mathbb{Y}^n$ expressed by $\mathcal{B}_{n,t}(f_n) \triangleq $
\begin{equation} \label{eq_lemma_achievability_1}
\left \{ (z,y_1^n) \in 
\{1,\dots,2^{nR} \} 
\times \mathbb{Y}^n:   \frac{P_{f_n(X_1^n)Y_1^n}(\doublevecn )}{Q_{f_n(X_1^n)Y_1^n}(\doublevecn)} > e^{nt}\right \}.
\end{equation}
$\mathcal{B}_{n,t}(f_n)$ is parametrized in terms of $t$, $n$ and $f_n$. Let us denote by $\phi_{n,t}(\cdot): \{1,\dots,2^{nR} \} \times \mathbb{Y}^n \mapsto \{0,1 \}$ the induced test (or decision rule) 
such that $\phi_{n,t}^{-1}(\left\{ 0 \right\})=\mathcal{B}_{n,t}(f_n)$. 
Then the \textsc{Type I} error probability for the pair $(f_n,\phi_{n,t})$ is given by
\begin{equation} \label{eq_lemma_achievability_2}
P_0(f_n,\phi_{n,t})= P_{f_n(X_1^n)Y_1^n} (\mathcal{B}^c_{n,t}(f_n)). 
\end{equation}
By construction of the pair $(f_n,\phi_{n,t})$, an upper bound for the \textsc{Type II} is obtained by
\begin{equation}
P_1(f_n,\phi_{n,t}) = Q_{f_n(X_1^n)Y_1^n}\left (	\mathcal{B}_{n,t}(f_n) \right )  \leq e^{-nt}.
\label{typeIIerrorcoundR}
\end{equation}
Then, for any finite $n>0$ and $\epsilon_n>0$,  finding an achievable \textsc{Type II} error exponent from this construction (and the bound in (\ref{typeIIerrorcoundR})) reduces to solve the following 
problem:
\begin{equation}
t^{\ast}_n(\epsilon_n)\triangleq \sup_{f_n \text{ encoder of rate }R} \sup_t \{t: 
P_{f_n(X_1^n)Y_1^n} (\mathcal{B}^c_{n,t}(f_n)) \leq \epsilon_n \}.
\label{optimumt}
\end{equation}

Note that  $f_n$ 
breaks the i.i.d. structure of the problem,  then determining $t^{\ast}_n(\epsilon_n)$ is not a simple task. 
We will derive a lower bound for $t^{\ast}_n(\epsilon_n)$ using a finite block analysis approach. For this, let us consider a fixed $l\geq 1$ and let us consider $\tilde{f}_l$ an encoder of length $l$, i.e. $\tilde{f}_l: \mathbb{X}^l \rightarrow \{1,\dots,2^{lR} \}$.  The idea is to decompose $X^n_1$ in segments of finite length to use the induced block  i.i.d. structure when $n$ tends to infinity. More precisely,  we construct an encoder that we denote by $\tilde{f}_{n,l}$ applying the function $\tilde{f}_l$ $k$-times to every sub-block of length $l$, assuming for the moment that $n=kl$, i.e., 
 \begin{align} \label{eq_lemma_achievability_3}
\tilde{f}_{n,l}(x_1,\dots,x_l,x_{l+1},\dots,x_{2l},\dots,x_{l(k-1)+1},\dots,x_{kl})\triangleq \nonumber \\
(\tilde{f}_l(x_1,\dots,x_l),\tilde{f}_l(x_{l+1},\dots,x_{2l}),\dots,\tilde{f}_l(x_{l(k-1)+1},\dots,x_{kl})).
 \end{align}
In the use of the set $\mathcal{B}_{n,t}(\tilde{f}_{n,l})$ in (\ref{eq_lemma_achievability_1}), 
it will be convenient to parametrize $t$  relative to the reference value $ \frac{1}{l}\mathcal{D}(P_{\tilde{f}_l(X_1^l)Y_1^l}\|Q_{\tilde{f}_l(X_1^l)Y_1^l})$ that is a function of $\tilde{f}_l$.
More precisely, let us define  
$$
t_{\delta}\triangleq   \frac{1}{l}\mathcal{D}(P_{\tilde{f}_l(X_1^l)Y_1^l}\|Q_{\tilde{f}_l(X_1^l)Y_1^l}) -\delta, 
$$ for any $\delta >0$.  Using the $l$-block structure of $\tilde{f}_{n,l}$, the \textsc{Type I} error in (\ref{eq_lemma_achievability_2}) of the pair $(\tilde{f}_{n,l}, \phi_{n,t_\delta})$ can be expressed by: 
\begin{equation}
P_{\tilde{f}_{n,l}(X_1^n)Y_1^n} \left (	\mathcal{B}^c_{n,t_{\delta}}(\tilde{f}_{n,l})\right ),
\end{equation}
where $\mathcal{B}^c_{n,t_{\delta}}(\tilde{f}_{n,l}))$ has the elements $z_1^k,y_1^n \in \{1,\dots,2^{lR} \}^k \times \mathbb{Y}^n$ satisfying that 
\begin{equation}
\Biggr \rvert \hat{\mathcal{D}}(P_{\tilde{f}_l(X_1^l)Y_1^l}\|Q_{\tilde{f}_l(X_1^l)Y_1^l})-\mathcal{D}(P_{\tilde{f}_l(X_1^l)Y_1^l}\|Q_{\tilde{f}_l(X_1^l)Y_1^l}) \Biggr \rvert \geq l\delta,  \label{eq_lemma_achievability_4}
\end{equation}
where $\hat{\mathcal{D}}(P_{\tilde{f}_l(X_1^l)Y_1^l}\|Q_{\tilde{f}_l(X_1^l)Y_1^l})\triangleq$
$$ \frac{1}{k}\sum_{i=1}^k \log  \left (\frac{P_{\tilde{f}_l(X_1^l)Y_1^l}(\{z_i,y_{k(i-1)+1}^{ki}\} )}{Q_{\tilde{f}_l(X_1^l)Y_1^l}(\{z_i,y_{k(i-1)+1}^{ki} \})} \right )
$$ 
denotes the empirical divergence. We will use a concentration inequality to bound the probability of the deviation event in (\ref{eq_lemma_achievability_4}). 
To this end,  let us introduce the notation:  $u_i=(z_i,y_{l(i-1)+1},\dots,y_{il})\in \{1,\dots,2^{lR} \}\times \mathbb{Y}^l$ and
\begin{equation} \label{eq_lemma_achievability_5}
 g(u_1,\dots,u_i,\dots,u_k)\triangleq \frac{1}{k} \sum_{j=1}^k \log \left (\frac{P_{\tilde{f}_l(X_1^l)Y_1^l}(\{u_j\} )}{Q_{\tilde{f}_l(X^l_1)Y_1^l}(\{u_j \})} \right ), 
\end{equation}
where it follows that for any $k>0$ and  $\forall i\in \{1,\dots,k \}$:
\begin{align} \label{eq_lemma_achievability_6}
&
\sup_{\begin{subarray}{c} u_1,\dots,u_i,\bar{u}_i,\dots,u_k\\
\in \tilde{f}_l(\mathbb{X}^l) \times \mathbb{Y}^l
\end{subarray}} \Biggr \rvert  g(u_1,\dots,u_i,\dots,u_k)
- g(u_1,\dots,\bar{u}_i,\dots,u_k) \Biggr \rvert  \nonumber \\ &  
\leq  \frac{2}{k} C(\tilde{f}_l,P_{XY}),
\end{align}
where $C(\tilde{f}_l,P_{XY})\triangleq \sup_{\doublevecl \in \tilde{f}_l(\mathbb{X}^l) \times \mathbb{Y}^l}\limits \Biggr \rvert \log\left( \frac{P_{\tilde{f}_l(X_1^l)Y_1^l}(\{ \doublevecl \})}{Q_{\tilde{f}_l(X_1^l)Y_1^l}(\{ \doublevecl \})}\right)\Biggr \rvert.$ From the bounded difference inequality \cite[Theorem 2.2]{devroye2012combinatorial}, we have that
\begin{equation} 
 P_{\tilde{f}_{n,l}(X_1^n)Y_1^n} \left (	\mathcal{B}^c_{n,t_{\delta}}(\tilde{f}_{n,l}) \right ) \leq \exp \left ( \frac{-k(l\delta)^2}{2C^2(\tilde{f}_l,P_{XY}) }\right ).
\label{boundep}
\end{equation}
Finally, from (\ref{optimumt}), a lower bound for $t^{\ast}_n(\epsilon_n)$ can be obtained from (\ref{boundep}) by making $\delta$ (that we denote by $\tilde{\delta}_{n,l}(\epsilon_n)$ in (\ref{eq_lemma_achievability_8})) the solution of the following condition:
\begin{equation} \label{eq_lemma_achievability_8}
\exp \left ( \frac{-k(l\tilde{\delta}_{n,l}(\epsilon_n))^2}{2 C^2(\tilde{f}_l,P_{XY})}\right )=\epsilon_n.
\end{equation}
Consequently, we have that 
\begin{align}
t^{\ast}_{n}(\epsilon_n)\geq & \underbrace{\frac{1}{l}\mathcal{D}(P_{\tilde{f}_l(X_1^l)Y_1^l}\|Q_{\tilde{f}_l(X_1^l)Y_1^l}) -\tilde{\delta}_{n,l}(\epsilon_n)}_{t_{\tilde{\delta}_{n,l}(\epsilon_n)}=}  
\label{boundupperR}
\end{align}
where from (\ref{eq_lemma_achievability_8}),
\begin{equation}
\tilde{\delta}_{n,l}(\epsilon_n)=\sqrt{\frac{2\ln(1/\epsilon_n)}{nl}}\cdot C(\tilde{f}_l,P_{XY}).
\label{deltaen}
\end{equation}
Finally, replacing the bound of $t^{\ast}_{n}(\epsilon_n)$  in (\ref{boundupperR}) at the exponential term in (\ref{typeIIerrorcoundR}) and taking logarithm, we have that:
\begin{dmath}
\xi(R)- \left( -\frac{1}{n}\log  P_1(\tilde{f}_{n,l},\phi_{n,t_{\tilde{\delta}_{n,l}(\epsilon_n)}}) \right )\leq \left [ \xi(R)-  \frac{1}{l}\mathcal{D}(P_{\tilde{f}_l(X_1^l)Y_1^l}\|Q_{\tilde{f}_l(X_1^l)Y_1^l})\right ] +\tilde{\delta}_{n,l}(\epsilon_n).
\label{principaldifference}
\end{dmath}

\begin{remark}
Looking at (\ref{principaldifference}) and using (2.6) and Theorem 3 in \cite{ahlswede1986hypothesis}, $\forall \gamma >0$, we can find a sufficient large $l^{\ast}$ and $f_l^{\ast}$ (function of $\gamma$) such that,
\begin{equation}
\xi(R)-\gamma <   \frac{\mathcal{D}(P_{\tilde{f}^{\ast}_l(X_1^{l^{\ast}})Y_1^{l^{\ast}}}\|Q_{\tilde{f}^{\ast}_l(X_1^{l^{\ast}})Y_1^{l^{\ast}}})}{l^{\ast}} < \xi(R).
\end{equation}
\end{remark}

Returning to the proof, we have that $\forall l > 0$, $\forall n > 0$ and  any $\epsilon_n > 0$
\begin{align}
& 
\xi(R)- \left( -\frac{1}{n}\log( \beta_n(\epsilon_n,R)) \right )  \nonumber \\
&\leq  \xi(R)- \left( -\frac{1}{n}\log( P_1(\tilde{f}_{n,l},\phi_{n,t_{\tilde{\delta}_{n,l}(\epsilon_n)}})) \right ) \nonumber \\
&\leq  \xi(R)-  \frac{1}{l}\mathcal{D}(P_{\tilde{f}_l(X_1^l)Y_1^l}\|Q_{\tilde{f}_l(X_1^l)Y_1^l}) +\tilde{\delta}_{n,l}(\epsilon_n) \nonumber \\
&= \left ( \max_{\begin{subarray}{c} U: U \mkv X \mkv Y \\ 
I(U;X) \leq R \\ |\mathbb{U}| \leq |\mathbb{X}|+1
\end{subarray}}I(U;Y)-  \frac{1}{l}I(\tilde{f}_l(X_1^l);Y_1^l)\right ) +\tilde{\delta}_{n,l}(\epsilon_n).
\label{overhead_achievable}
\end{align} 
The first inequality is from the fact that $\beta_n(\epsilon_n,R))\leq P_1(\tilde{f}_{n,l},\phi_{n,t_{\tilde{\delta}_{n,l}(\epsilon_n)}})$, the second from (\ref{principaldifference}), and  
the last equality from the definition of $\xi(R)$ in Lemma \ref{th_ahlswede&csiszarlimit}, expressing the divergence as a mutual information \cite{cover2012elements}. 

It is worth noting that the bound in (\ref{overhead_achievable}) is valid for an arbitrary $l>0$. Considering that we know an expression for $\tilde{\delta}_{n,l}(\epsilon_n)$ from (\ref{deltaen}), we can address the problem of  finding the best upper bound, i.e., the $l$ that offers the best compromise between 
the two terms in the RHS of (\ref{overhead_achievable}). For that, we need to focus on: 
\begin{equation}
\max_{\begin{subarray}{c} U: U \mkv X \mkv Y \\ 
I(U;X) \leq R \ |\mathbb{U}| \leq |\mathbb{X}|+1
\end{subarray}}I(U;Y)-  \max_{\tilde{f}_l: \mathbb{X}^l \rightarrow \{1,...,2^{lR} \}}\frac{1}{l}I(\tilde{f}_l(X_1^l);Y_1^l), 
\label{bottleneck}
\end{equation}
which corresponds to the non-asymptotic analysis of the information bottleneck problem \cite{tishby2000information}. 
This coding problem can be viewed as a classical rate-distortion (fixed-rate) lossy source coding problem with the log-loss as the distortion function \cite{courtade2014multiterminal}. More precisely, (\ref{bottleneck}) can be expressed by:
\begin{equation}
 \min_{\tilde{f}_l: \mathbb{X}_1^l \rightarrow \{1,\dots,2^{lR} \}}\frac{1}{l}H(Y_1^l|\tilde{f}_l(X_1^l))-\min_{\begin{subarray}{c} U: U \mkv X \mkv Y  \\ 
I(U;X)\leq R \ |\mathbb{U}| \leq |\mathbb{X}|+1
\end{subarray}}H(Y|U).
\label{bounddistortion}
\end{equation}

The following Lemma connects the expression in (\ref{bounddistortion}) with an instance of the classical rate distortion problem \cite{berger1971rate}.
\begin{lemma}
\label{zhanglogloss}
\begin{equation}
\begin{split}
\frac{1}{l}H(Y_1^l|\tilde{f}_l(X_1^l) & \leq  D(R)-\frac{\partial}{\partial R}D(R)\frac{\log (l)}{2l} + o\left( \frac{\log l}{l}\right ),
\label{overhead_zhang}
\end{split}
\end{equation}
where $D(R)$ is the noisy distortion-rate function given by
\begin{equation}
D(R)= \min_{\begin{subarray}{c} U: U \mkv X \mkv Y  \\ 
I(U;X)\leq R \ |\mathbb{U}| \leq |\mathbb{X}|+1
\end{subarray}}\limits H(Y|U).
\label{noisyrate}
\end{equation}
\end{lemma}
The proof is presented in Appendix \ref{appendC}. Consequently, from (\ref{overhead_zhang}) we have that the expression in (\ref{bottleneck}) is upper bounded by
$-\frac{\partial}{\partial R}D(R)\frac{\log (l)}{2l} + o\left( \frac{\log (l)}{l}\right ).$
Applying this result to (\ref{overhead_achievable}), it follows that
\begin{align}
 & \xi(R)- \left( -\frac{1}{n}\log( \beta_n(\epsilon_n,R)) \right ) \nonumber \\
 & \leq  -\frac{\partial}{\partial R}D(R)\frac{\log (l)}{2l} +\tilde{\delta}_{n,l}(\epsilon_n)+ o\left( \frac{\log (l)}{l}\right ).
\label{overhead_final}
\end{align} 
To obtain a more explicit dependency of $\tilde{\delta}_{n,l}(\epsilon_n)$ on $l$ we use the following result:
\begin{proposition}
\label{Scaling}
Let us consider two arbitrary probability distributions $\mu, \rho \in \mathbb{P}(\mathbb{X})$, an arbitrary encoder $f_n:\mathbb{X}\rightarrow \{1,\dots,n\}$. and its induced partition of $\mathbb{X}$ given by
$\pi_n=\{A_{i,n} \triangleq f_n^{-1}(\{i\}): i \in \{1,\dots,n \} \}$,
then
\begin{equation}
\sup_{A\in \pi_n}\frac{\mu(A)}{\rho(A)}\leq \sup_{x\in \mathbb{X}}\frac{\mu(\{x\})}{\rho(\{x\})}.
\end{equation}
\end{proposition}
The proof is presented in Appendix \ref{appendA}.
%

From Proposition \ref{Scaling}, we obtain that: 
\begin{align}
\tilde{\delta}_{n,l}(\epsilon_n)& =\sqrt{\frac{2\ln(1/\epsilon_n)}{nl}}\cdot C(\tilde{f}_l,P_{XY}) \nonumber \\
 & \leq \sqrt{\frac{2l\ln(1/\epsilon_n)}{n}}\cdot C(P_{XY}).
 \label{newdelta}
\end{align}
Using (\ref{newdelta}), 
the problem reduces to minimize the RHS of (\ref{overhead_final}) as long as $(\epsilon_n)_n$ tends to zero at a sub-exponential rate, for which the assumption that $\left ( \frac{1}{\epsilon_n} \right )_n$ is $o(e^{rn})$ for any $r>0$ is central. In fact, it is sufficient to consider any sequence $(l_n)_n$ of integers such that $(1/l_n)_n$ is $o(1)$ and $(l_n)_n$ is $o\left ( \frac{n}{\ln (1/\epsilon_n)}\right)$, from which we conclude that  
\begin{equation}
\lim_{n\rightarrow \infty}\inf - \frac{1}{n}\log(\beta_n(\epsilon_n,R))\geq \xi(R).
\end{equation}

\end{IEEEproof}

Conversely,  we have the following result:
\begin{theorem}\label{lemmadehan}
Let us assume that $\epsilon_n>0$ for all $n$ and that  $(1/{\epsilon_n})_n=o(e^{rn})$ for any $r >0$, then
\begin{equation}
\lim_{n\rightarrow \infty}\sup - \frac{1}{n}\log(\beta_n(\epsilon_n,R))\leq \xi(R).
\end{equation}
\end{theorem}

\begin{IEEEproof}
Let us consider a fixed-rate encoder $f_n: \mathbb{X}^n \rightarrow \{1,\dots,2^{nR}\}$ of rate $R$.
We begin by using \cite[Lemma 4.1.2]{koga2013information}, which states that for all $t>0$ and $\forall \mathcal{A}_n \subset f_n(\mathbb{X}^n) \times \mathbb{Y}^n$
\begin{dmath}
P_{f_n(X_1^n)Y_1^n}(\mathcal{A}^c_n)+e^{nt}Q_{f_n(X_1^n)Y_1^n}(\mathcal{A}_n)\geq P_{f_n(X_1^n)Y_1^n}\left( \mathcal{B}_{n,t}^c(f_n)\right ),
\label{hanlemma}
\end{dmath}
where as before $\mathcal{B}_{n,t}(f_n)=$
$$ \left \{(z,y_1^n) \in f_n(\mathbb{X}^n) \times \mathbb{Y}^n : \frac{P_{f_n(X_1^n)Y_1^n}(\{\doublevecn \})}{Q_{f_n(X_1^n)Y_1^n}(\{\doublevecn\})} > e^{nt}\right \}.$$
(\ref{hanlemma}) 
is valid for any 
decision rule acting on $(f_n(X_1^n),Y_1^n)$.  The rest of the argument focuses on finding a lower bound to the RHS of (\ref{hanlemma}). The latter can be done by considering the following 
function $i(x_1^n,y_1^n)=\log \left( \frac{P_{Y_1^n|f_n(X_1^n)(y_1^n|f_n(x^n_1))}}{P_{Y_1^n}(y^n_1)}\right)$ and the fact that $\forall q \geq 1$
\begin{equation}
\mathbb{E}_{(X^n_1,Y^n_1)\sim P^n_{XY}}(i(X_1^n,Y_1^n)^q) \leq q!(4n^2C(P_{XY})^2)^q.
\end{equation}
Using $\mathcal{B}_{n,t}(f_n)$, it is useful to write $t =\frac{I(f_n(X_1^n);Y_1^n)}{n}+s$,  
then $P_{f_n(X_1^n)Y_1^n}\left( \mathcal{B}_{n,t}^c(f_n)\right )=$
\begin{equation}\label{eq_unlabed_1}
P^n_{XY} \left( \left\{(x^n_1,y^n_1): i(x_1^n,y_1^n)-\mathbb{E}(i(X_1^n,Y_1^n))\leq ns \right\}  \right), 
\end{equation}
where the expected values in (\ref{eq_unlabed_1}) assumes that $(X^n_1,Y^n_1)\sim P^n_{XY}$.
Using the bound on the variance of $i(X_1^n,Y_1^n)$, we can use the moment concentration inequality \cite[Theorem 2.1]{boucheron2013concentration} to obtain: 
\begin{align}
& P^n_{XY} \left(  \left\{ i(x_1^n,y_1^n)-\mathbb{E}(i(X_1^n,Y_1^n))\leq ns \right\}  \right) \\
 & \geq 1-e^{-s^2/(32C(P_{XY})^2)} \nonumber.
\end{align}
Combining this with (\ref{hanlemma}), it follows that for any $s>0$ and any set $\mathcal{A}_n \subset f_n(\mathbb{X}^n) \times \mathbb{Y}^n$
\begin{align}
& P_{f_n(X_1^n)Y_1^n}(\mathcal{A}^c_n)+e^{n\left (\frac{I(f_n(X_1^n);Y_1^n)}{n}+s \right )}Q_{f_n(X_1^n)Y_1^n}(\mathcal{A}_n) \nonumber \\
& \geq 1-e^{-s^2/(32C(P_{XY})^2)}.
\label{hanlemmanew}
\end{align}

At this point, we introduce the restriction on the \textsc{Type I} error in the analysis. Let us consider an arbitrary $\mathcal{A}_n$ such that $P_{f_n(X_1^n)Y_1^n}(\mathcal{A}^c_n) \leq \epsilon_n$. Then we have that:
\begin{align}
&e^{n\left (\frac{I(f_n(X_1^n);Y_1^n)}{n}+s \right )}Q_{f_n(X_1^n)Y_1^n}(\mathcal{A}_n) \nonumber \\
&\geq 1-e^{-s^2/(32C(P_{XY})^2)}-\epsilon_n.
\label{hanlemmanew2}
\end{align}
Taking logarithm at both sides of (\ref{hanlemmanew2}) for any $s$ satisfying the admisible condition 
$\epsilon_n < 1 -e^{-\frac{s^2}{32C(P_{XY})^2}}$, it follows that
\begin{align}
& \frac{I(f_n(X_1^n);Y_1^n)}{n}- \left( -\frac{1}{n}\log(Q_{f(X_1^n)Y_1^n}(\mathcal{A}_n)) \right )\nonumber \\
& \geq -s+ \frac{\log\left ( 1-\epsilon_n -e^{-\frac{s^2}{32C(P_{XY})^2}}\right )}{n}.
\label{finalbound}
\end{align}
Since both $f_n$ and the set $\mathcal{A}_n$ are arbitrary in (\ref{finalbound}), the bound is valid for the pair $(f^{\ast}_n,\phi^{\ast}_n)$ such that $Q_{f_n^{\ast}(X_1^n)Y_1^n}(\mathcal{A}^{\ast}_n)=\beta_n(\epsilon_n,R)$.  In addition $\frac{I(f_n(X_1^n);Y_1^n)}{n} \leq \xi(R)$ by definition (see (2.5) in \cite{ahlswede1986hypothesis}), then for all $s>4C(P_{XY})\sqrt{2\ln (1/1-\epsilon_n)}$ it follows that
\begin{align}\label{eq_final_equation}
\xi(R)+&\frac{1}{n}\log( \beta_n(\epsilon_n,R))\geq\nonumber\\ 
&-s+ \frac{\log\left ( 1-\epsilon_n -e^{-\frac{s^2}{32C(P_{XY})^2}}\right )}{n}.
\end{align}
At this point, we use the assumption that $\lim_{n \rightarrow \infty} \epsilon_n = 0$, which implies that 
there is a sequence $(s_n)_n$ that is $\mathcal{O}(\sqrt{\log(n)/n})$ for which (\ref{eq_final_equation}) evaluated at $s=s_n$ holds
for any $n$, which implies that
\begin{equation}
\lim_{n\rightarrow \infty}\sup - \frac{1}{n}\log(\beta_n(\epsilon_n,R))\leq \xi(R).
\end{equation}
\end{IEEEproof}

\subsection{Proof of Theorem \ref{theorem2}}
\label{proofresult2}
\begin{IEEEproof}
The proof can be divided in two independent parts from the analysis obtained in Theorems \ref{lemma_achievability} and \ref{lemmadehan}. On the one hand, we have an upper bound obtained by optimizing the RHS of (\ref{overhead_final}) with respect to the blocklength $l$. More precisely, we have the following inequality:
\begin{align}
 &\xi(R)- \left( -\frac{1}{n}\log( \beta_n(\epsilon_n,R)) \right )\nonumber \\
 &\leq  -\frac{\partial}{\partial R}D(R)\frac{\log l}{2l} +\sqrt{\frac{2l\ln(1/\epsilon_n)}{n}} C(P_{XY}) + o\left( \frac{\log l}{l}\right ),
\label{match}
\end{align}
where $C(P_{XY})\triangleq \sup_{(x,y) \in \mathbb{X} \times \mathbb{Y}}\Biggr \rvert \log\left( \frac{P_{XY}(\{ (x,y) \})}{Q_{XY}(\{ (x,y) \})}\right)\Biggr \rvert$.  This expression depends on $\epsilon_n$ and it is valid for all $l\geq 1$. Then the tighest bound from (\ref{match}), reduces to find $l_n^{\ast}$ solution of:
\begin{equation}
\frac{\log l_n^{\ast}}{l_n^{\ast}}
\approx
\sqrt{\frac{l_n^{\ast}\ln(1/\epsilon_n)}{n}}. 
\label{matchalpha}
\end{equation}
To address this problem, we consider $l_n=n^{\alpha}$ to look for this optimal $\alpha$ (function of $\epsilon_n$). This is the consequence of assuming that the condition in (\ref{matchalpha}) holds, which reduces to:
\begin{equation}
\frac{\log n^{\alpha}}{n^{\alpha}} 
\approx 
\sqrt{\frac{n^{\alpha}\ln(1/\epsilon_n)}{n}}.
\label{match2}
\end{equation}
To solve (\ref{match2}), we move into the specific cases for ($\epsilon_n$) stated in Theorem \ref{theorem2}. We have three different scenarios:\\
a) $(\epsilon_n)_n=(1/n^p)_n$ with $p> 0$: The condition (\ref{match2}) reduces to
\begin{equation}
\frac{\alpha\log n}{n^{\alpha}}
\approx 
\sqrt{\frac{n^{\alpha}p\ln(n)}{n}}, 
\label{match3}
\end{equation}
where (non considering the logarithmic term) the equilibrium is obtained with $\alpha^{\ast}=1/3$, which makes the upper bound in (\ref{match}) of the form:
\begin{align}
 & \xi(R)- \left( -\frac{1}{n}\log( \beta_n(\epsilon_n,R)) \right ) \nonumber \\ & \leq  -\frac{\partial D(R)}{\partial R}\frac{\log n}{6n^{1/3}} +\sqrt{\frac{2p\ln(n)}{n^{2/3}}}  C(P_{XY}) + o\left( \frac{\log n}{n^{1/3}}\right )\nonumber \\
 & =  \left [-\frac{\partial D(R)}{\partial R}\cdot \frac{1}{6} + o\left(1\right )\right ]\left (\frac{\log n}{n^{1/3}} \right ). 
 \label{solution}
\end{align}

b) $(\epsilon_n)_n=(1/e^{n^p})_n$ with $p\in( 0,1)$: Following the previous approach,  we solve
\begin{equation}
 \frac{\alpha\log n}{n^{\alpha}}
 \approx 
 \sqrt{\frac{n^{\alpha}n^p}{n}}, 
\label{match4}
\end{equation}
resulting in $\alpha^{\ast}=(1-p)/3$. This choice offers the  bound
\begin{align}
 & \xi(R)- \left( -\frac{1}{n}\log( \beta_n(\epsilon_n,R)) \right ) \nonumber \\ & \leq  -\frac{\partial D(R)}{\partial R}\frac{(1-p)\log n}{6n^{(1-p)/3}} +\frac{\sqrt{2} C(P_{XY})}{n^{(1-p)/3}}+ o\left( \frac{\log n}{n^{(1-p)/3}}\right )\nonumber \\
 & = \left [-\frac{\partial D(R)}{\partial R}\frac{(1-p)}{6} + o\left( 1\right )\right ] \left (\frac{\log n}{n^{(1-p)/3}} \right ).
\end{align}

c) $(\epsilon_n)_n=(1/\log(n))_n$: The  matching condition reduces to find $\alpha$ such that
\begin{equation}
\frac{\alpha\log n}{n^{\alpha}} 
\approx 
\sqrt{\frac{n^{\alpha}\ln(\log(n))}{n}}. 
\label{match5}
\end{equation}
It is simple to show that, as in the polynomial regime, the approximated solution is $\alpha^{\ast}=1/3$, which offers the following upper bound:
\begin{align}
 & \xi(R)- \left( -\frac{1}{n}\log( \beta_n(\epsilon_n,R)) \right ) \nonumber \\ & \leq  -\frac{\partial D(R)}{\partial R}\frac{\log n}{6n^{1/3}} +\sqrt{\frac{2\ln(\log(n))}{n^{2/3}}} C(P_{XY})+ o\left( \frac{\log n}{n^{1/3}}\right )\nonumber \\
 & =  \left [-\frac{\partial D(R)}{\partial R}\cdot\frac{1}{6} + o\left(1\right )\right ]\left (\frac{\log n}{n^{1/3}} \right ). 
 \label{solution3}
\end{align}

For the lower bound, we use the following inequality from the proof of Theorem \ref{lemmadehan} (see (\ref{eq_final_equation})):
\begin{align}
\xi(R)&- \left( -\frac{1}{n}\log( \beta_n(\epsilon_n,R)) \right )\geq \nonumber\\ 
&-s+ \frac{\log\left ( 1-\epsilon_n -e^{-\frac{s^2}{32C(P_{XY})^2}}\right )}{n}.
\label{finaloverheadconverse}
\end{align}
 This inequality is valid for any $s\in \R$ such that $ 1-\epsilon_n -e^{-\frac{s^2}{32C(P_{XY})^2}}>0$ or, equivalently, for all $s$ such that $s>4C(P_{XY})\sqrt{2\ln (1/1-\epsilon_n)}$. At this point, it is convenient to define $h_n(s)\triangleq 1-\epsilon_n -e^{-\frac{s^2}{32C(P_{XY})^2}}$ in the domain $s>4C(P_{XY})\sqrt{2\ln (1/1-\epsilon_n)}$. Then (\ref{finaloverheadconverse}) can be expressed in terms of $h_n(s)$ by
\begin{align}
& \xi(R)- \left( -\frac{1}{n}\log( \beta_n(\epsilon_n,R)) \right )\nonumber \\
& \geq -4C(P_{XY})\sqrt{2\ln\left( \frac{1}{1-\epsilon_n-h_n(s)}\right)}-\frac{\log(1/h_n(s))}{n},
\label{conditiontwo}
\end{align}
where $h_n(s) >0$ if $s>4C(P_{XY})\sqrt{2\ln (1/1-\epsilon_n)}$. We notice that as $(\epsilon_n)_n$ is $o(1)$ (function of $n$) the first term on the RHS of (\ref{conditiontwo}) tends to zero if, and only if,  $(h_n(s))_n$ is $o(1)$. On the other hand, $(\log(1/h_n(s)))_n$ needs to be $o(n)$ to make the second terms on the RHS of (\ref{conditiontwo}) vanishing to zero with $n$. Then, there is a regime on the asymptotic behavior of $(h_n(s))_n$ where the bound in (\ref{conditiontwo}) is meaningful.  

More precisely, for any finite $n$, we will address the problem of finding $s \in (4C(P_{XY})\sqrt{2\ln (1/1-\epsilon_n)},\infty)$,  or equivalently finding $h_n(s) \in (0,1)$,  that offers the best lower bound from (\ref{conditiontwo}). On the specifics, as $(\epsilon_n)_n$ and $(h_n(s))_n$ go to zero with $n$, for the first term $-4C(P_{XY})\sqrt{2\ln\left( \frac{1}{1-\epsilon_n-h_n(s)}\right)}$
 a Taylor expansion around 1 is used to aproximate the function. In particular, it follows that:
\begin{align}
& -4\sqrt{2}C(P_{XY})\sqrt{\ln\left( \frac{1}{1-\epsilon_n-h_n(s)}\right)} \nonumber \\
& \geq -2\sqrt{2}C(P_{XY})\sqrt{\epsilon_n+h_n(s)}\frac{\sqrt{4-5 (\epsilon_n+h_n(s))}}{1-\epsilon_n-h_n(s)}\nonumber\\
& \geq  -2\sqrt{2}C(P_{XY})\sqrt{\epsilon_n+h_n(s)}\frac{\sqrt{4}}{1/2}\nonumber\\ 
&= -8\sqrt{2}C(P_{XY})\sqrt{\epsilon_n+h_n(s)}, 
\label{conditionthree}
\end{align}
where the last inequality is obtained eventually as $(\epsilon_n+h_n(s))_n$ is $o(1)$. Then, from (\ref{conditiontwo}) and (\ref{conditionthree}), the optimal lower bound reduces to find the optimal balance between $-8\sqrt{2}C(P_{XY})\sqrt{\epsilon_n+h_n(s)}$ and $\frac{\log(1/h_n(s))}{n}$.  It is important to note that $-8\sqrt{2}C(P_{XY})\sqrt{\epsilon_n+h_n(s)}$ tends to zero at a velocity that is proportional to how fast $(h_n(s))_n$ tends to zero, as long as, $(h_n(s))_n$ is $o(\epsilon_n)$, otherwise,  the velocity is dominated by $\mathcal{O}(\sqrt{\epsilon_n})$, which is independent of $(h_n(s))_n$. 
On the other hand, the second term $(\log(1/h_n(s)))_n$ tends to zero at a rate that is inversely proportional to the velocity at which 
$(h_n(s))_n$ goes to zero.  Therefore, the balance is function of $(\epsilon_n)_n$.
We recognize two regimes for this optimization problem:
\begin{itemize}
\item[1-] If for some $K>0$ we have that $\sqrt{2\epsilon_n}\geq K \frac{\log(1/\epsilon_n)}{n}$, eventually in $n$, then the solution of the optimization problem is achieved when $(h_n(s))_n \approx (\epsilon_n)_n$ ({\bf Regime 1});
\item[2-] Otherwise, if $(\sqrt{2\epsilon_n})_n$ is  $o\left (\frac{\log(1/\epsilon_n)}{n} \right )$, then the solution of the optimization problem implies that $(\epsilon_n)_n$ is $o(h_n(s))$ ({\bf Regime 2}).
\end{itemize}
Finally, to obtain the upper bound, we need to evaluate $(\epsilon_n)_n$ in the different scenarios stated in Theorem \ref{theorem2}. 
\begin{itemize}
\item $(\epsilon_n)_n=(1/\log(n))_n$: Regime 1 is met, then we choose $h_n(s)=\epsilon_n$. This implies that 
\begin{align}
&\xi(R)- \left( -\frac{1}{n}\log( \beta_n(\epsilon_n,R)) \right )\nonumber \\ & \geq \frac{-16C(P_{XY})}{\sqrt{\log(n)}}-\frac{\log(\log(n))}{n} \nonumber \\
& = \left (-16C(P_{XY})-\frac{\log(\log(n))\sqrt{\log(n)}}{n} \right )\frac{1}{\sqrt{\log(n)}}\nonumber\\
& = \left (-16C(P_{XY})-o(1) \right )\left (\frac{1}{\sqrt{\log(n)}}\right ).
\label{match9}
\end{align}

\item $(\epsilon_n)_n=(1/n^p)_n$ with $2>p> 0$: Regime 1 is met, then we choose $h_n(s)=\epsilon_n$. This implies that 
\begin{align}
&\xi(R)- \left( -\frac{1}{n}\log( \beta_n(\epsilon_n,R)) \right ) \nonumber \\ & \geq \frac{-16C(P_{XY})}{n^{p/2}}-\frac{p\log(n)}{n} \nonumber \\
& = \left (-16C(P_{XY})-\frac{p\log(n)}{n^{1-p/2}} \right )\left (\frac{1}{n^{p/2}} \right ) \nonumber\\
& = \left (-16C(P_{XY})-o(1) \right )\left (\frac{1}{n^{p/2}}  \right ).
\label{match6}
\end{align}

\item $(\epsilon_n)_n=(1/n^p)_n$ with $p\geq 2$: Regime 2 is met, then we have to solve the following matching condition 
\begin{equation}
 \sqrt{\epsilon_n+h_n(s)} \approx \frac{\log(1/h_n(s))}{n}.
\label{match10}
\end{equation}
Assuming $h_n(s)=1/n^{\alpha}$, $\alpha \in (0,2]$, the equilibrium is obtained with $\alpha^{\ast}=2$. This implies that 
\begin{align}
&\xi(R)- \left( -\frac{1}{n}\log( \beta_n(\epsilon_n,R)) \right )\nonumber \\ & \geq -8\sqrt{2}C(P_{XY})\sqrt{n^{-p}+n^{-2}}-\frac{2\log(n)}{n} \nonumber \\
& = \left (-8\sqrt{2}C(P_{XY})\frac{\sqrt{n^{2-p}+1}}{\log(n)}-2\right )\left ( \frac{\log(n)}{n} \right ) \nonumber \\
& = \left (-o(1)-2\right )\left ( \frac{\log(n)}{n} \right ).
\label{match7}
\end{align}

\item $(\epsilon_n)_n=(1/e^{n^p})_n$ with $p\in( 0,1)$: Regime 2 is met, then we follow the same condition in (\ref{match10}). The equilibrium is obtained with $\alpha^{\ast}=2$. This implies that
\begin{align}
& \xi(R)- \left( -\frac{1}{n}\log( \beta_n(\epsilon_n,R)) \right )\nonumber \\ & \geq -8\sqrt{2}C(P_{XY})\sqrt{e^{-n^{p}}+n^{-2}}-\frac{\log(n)}{n} \nonumber \\
& = \left (-8\sqrt{2}C(P_{XY})\frac{\sqrt{e^{-n^{p}}n^2+1}}{\log(n)}-2\right )\left ( \frac{\log(n)}{n} \right ) \nonumber \\
& = \left (-o(1)-2\right )\left ( \frac{\log(n)}{n} \right ).
\label{match8}
\end{align}
\end{itemize}
\end{IEEEproof}

\subsection{Proof of Lemma \ref{zhanglogloss}}
\begin{IEEEproof}
\label{appendC}
Let us consider a family of probability distributions $P_{\lambda} \in \mathcal{P}(\mathbb{Y})$ indexed with a parameter $\lambda \in \Lambda$, where $\Lambda$ is some parametric space.  Given a vector of parameters $\lambda_1^n \in \Lambda^n$, the product probability distribution in $\mathcal{P}(\mathbb{Y}^n)$ is defined as
\begin{equation}
P_{\lambda^n}(\{y_1^n \})\triangleq \prod_{i=1}^nP_{\lambda_i}(\{y_i\}).
\end{equation}
Let $\rho(\lambda_1^n,Y_1^n):\Lambda^n \times \mathbb{Y}_1^n \rightarrow \R^{+} \cup \{ 0\} $ denote the logarithmic loss distortion defined by:
\begin{equation*}
\rho(\lambda_1^n,y_1^n)\triangleq -\frac{1}{n}\log P_{\lambda_1^n}(\{y_1^n \})=\sum_{i=1}^n-\frac{1}{n}\log P_{\lambda_i}(\{y_i\}).
\label{log-loss}
\end{equation*}
By construction $\rho(\lambda_1^n,y_1^n)$ is additive $(\rho(\lambda_1^n,y_1^n)= \sum_{i=1}^n\rho(\lambda_i,y_i))$ and then the following result holds: 
\begin{lemma}
\label{equivalent_log}
\cite[Lemma 1]{courtade2014multiterminal} Let $X_1^l,Y_1^l$ be a random vector with known joint distribution. For any function $\tilde{f}_l:\mathbb{X}^l \rightarrow \{1,...,2^{lR} \}$ and function $g:\{1,...,2^{lR} \} \rightarrow \Lambda^n$ such that $g(\tilde{f}_l(X_1^l))=\lambda_1^l$ it follows that
\begin{equation}
\mathbb{E}[\rho(g(u),Y_1^l)|\tilde{f}_l(X_1^l)=u] \geq \frac{1}{l}H(Y_1^l|\tilde{f}_l(X_1^l)=u).
\label{loginequality}
\end{equation}
Taking expectation on the two sides of (\ref{loginequality}) with repect to $X^l_1$, we get that
\begin{equation}
\mathbb{E}[\rho(g(\tilde{f}_l(X_1^l)),Y_1^l)] \geq \frac{1}{l}H(Y_1^l|\tilde{f}_l(X_1^l)).
\label{loginequality2}
\end{equation}
\end{lemma}
\begin{remark}
We observe that, if we identify the $\tilde{f}_l$ as an encoder and $g$ as the decoder, the term in the LHS of (\ref{loginequality2}) corresponds to the noisy rate distortion function under the logarithmic loss. Then, for the purpose of the following result, it is convenient to redefine the distortion function $\tilde{\rho}(x_1^l,\lambda_1^l): \mathbb{X}_1^l \times \Lambda_1^l\rightarrow \R \cup \{0 \}$ as
\begin{equation}
\tilde{\rho}(x_1^l,\lambda_1^l)\triangleq \mathbb{E}[\rho(\lambda_1^l,Y_1^l)|X_1^l=x_1^l].
\end{equation}
Denoting $\lambda_i=g_i(\tilde{f}_l(x_1^l))$  and $g_i$ is the $i$th component of $g$, we observe that $\tilde{\rho}(x_1^l,\lambda_1^l)= \sum_{i=1}^l \tilde{\rho}(x_i,\lambda_i) $ is additive. 
\end{remark}
Finally, using the previous observation, we can use $\tilde{f}_l$ as the encoder and $g_i$ as the decoder, to recover an instance of the rate distortion problem \cite{berger1971rate}. Therefore, from \cite[Theorem 3]{zhang1997redundancy}, we obtain that
\begin{equation*}
\begin{split}
\frac{1}{l}H(Y_1^l|\tilde{f}_l(X_1^l))&\leq \mathbb{E}_{X\sim P^l_{X}}[\tilde{\rho}(X_1^l,\lambda_1^l)]\\
& \leq  D(R)-\frac{\partial}{\partial R}D(R)\frac{\log (l)}{2l} + o\left( \frac{\log (l)}{l}\right ),
\label{overhead_zhang2}
\end{split}
\end{equation*}
which concludes the result.
\end{IEEEproof}
\subsection{Finite-length Result for the Unconstrained Case}
\label{app_finite-length_unconstrained}
\begin{lemma} \label{th_strassen2009asymptotic} \cite{strassen2009asymptotic} 
Let us consider $\epsilon \in (0,1)$, then eventually in $n$ it follows that $-\frac{\log(\beta_n(\epsilon))}{n}=$
\begin{dmath*}
\mathcal{D}(P\|Q)+\sqrt{\frac{V(P\|Q)}{n}}\Phi^{-1}(\epsilon)+\frac{\log n}{2n} + \mathcal{O}\left (\frac{1}{n} \right),
\end{dmath*}
where $V(P\|Q)=\sum_{x \in \mathbb{X}}\limits P(\{x\})\left [\log \left (\frac{P(\{x\})}{Q(\{x\})} \right ) -\mathcal{D}(P\|Q) \right ]^2$.
\end{lemma}

A 
direct corollary of this result shows that  
$\left \rvert \mathcal{D}(P\|Q)- \left( -\frac{1}{n}\log(\beta_n(\epsilon)) \right )\right \rvert$ is  $\mathcal{O}\left(\frac{1}{\sqrt{n}}\right)$.

\subsection{Proof of Proposition \ref{Scaling}}
\begin{IEEEproof}
\label{appendA}
Given $A \in \pi_n$, we note that
\begin{equation} \label{eq_appendA_1}
\frac{\mu(\mathcal{A})}{\rho(\mathcal{A})}=\frac{\sum_{j=1}^{|A|}\limits \mu(\{j: j \in \mathcal{A} \})}{\sum_{j=1}^{|\mathcal{A}|}\limits \rho(\{j: j \in \mathcal{A} \})}.
\end{equation}
Then, given a collection of positive numbers $\{a_i: i \in \{1,\dots,n \} \}$ and $\{b_i: i \in \{1,\dots,n \} \}$, we use the following basic inequality
\begin{equation}
{\sum_{i=1}^n\limits a_i} \leq \max_i \left \{\frac{a_i}{b_i} \right \} {\sum_{i=1}^n\limits b_i}.
\label{calculus}
\end{equation}
Finally, since $\mathcal{A}$ is arbitrary in (\ref{eq_appendA_1}) and using the positiveness of the probability, we conclude the desired result.
\end{IEEEproof}

\ifCLASSOPTIONcaptionsoff
  \newpage
\fi

\bibliographystyle{IEEEtran}
\bibliography{referencias}

\begin{thebibliography}{10}
\providecommand{\url}[1]{#1}
\csname url@rmstyle\endcsname
\providecommand{\newblock}{\relax}
\providecommand{\bibinfo}[2]{#2}
\providecommand\BIBentrySTDinterwordspacing{\spaceskip=0pt\relax}
\providecommand\BIBentryALTinterwordstretchfactor{4}
\providecommand\BIBentryALTinterwordspacing{\spaceskip=\fontdimen2\font plus
\BIBentryALTinterwordstretchfactor\fontdimen3\font minus
  \fontdimen4\font\relax}
\providecommand\BIBforeignlanguage[2]{{%
\expandafter\ifx\csname l@#1\endcsname\relax
\typeout{** WARNING: IEEEtran.bst: No hyphenation pattern has been}%
\typeout{** loaded for the language `#1'. Using the pattern for}%
\typeout{** the default language instead.}%
\else
\language=\csname l@#1\endcsname
\fi
#2}}

\bibitem{Mahler_2007}
R.~Mahler, \emph{Statistical Multisources-Multitarget Information
  Fusion}.\hskip 1em plus 0.5em minus 0.4em\relax Norwood, MA, USA, 2007.

\bibitem{marano_2019}
S.~Marano and P.~K. Willet, ``Algorithm and fundamental limits for unlabeled
  detection using types,'' \emph{IEEE Transactions on Signal Processing},
  vol.~67, no.~8, pp. 2022--2035, 2019.

\bibitem{Wang_2018}
G.~Wang, J.~Zhu, R.~Blum, P.~K. Willet, S.~Marano, V.~Matta, and P.~Braca,
  ``Signal amplitude estimation and detection from unlabeled binary quantized
  samples,'' \emph{IEEE Transactions on Signal Processing}, vol.~66, no.~16,
  pp. 4291--4303, August 2018.

\bibitem{Zhu_2017}
J.~Zhu, H.~Cao, C.~Song, and Z.~Xu, ``Parameter estimation via unlabelled
  sensing using distributed sensors,'' \emph{IEEE Commun. Letter}, vol.~21,
  no.~10, pp. 2130--2133, 2017.

\bibitem{marano_2016}
S.~Marano and P.~K. Willet, ``The importance of being earnest: Social network
  with unknown agent quality,'' \emph{IEEE Transactions on Signal Processing},
  vol.~2, no.~3, pp. 306--320, September 2016.

\bibitem{unnikrishnam_2018}
J.~Unnikrishnam, S.~Haghighasthoar, and M.~Vetterli, ``Unlabeled sensing with
  random linear measurements,'' \emph{IEEE Transactions on Information Theory},
  vol.~64, no.~5, pp. 3237--3253, 2018.

\bibitem{liu_2018}
Z.~Liu and J.~Zhu, ``Signal detection from unlabeled ordered samples,''
  \emph{IEEE Commun. Letter}, vol.~22, no.~12, pp. 2431--2434, December 2018.

\bibitem{Haghighasthoar_2018}
S.~Haghighasthoar and G.~Caire, ``Signal recovery from unlabeled samples,''
  \emph{IEEE Transactions on Signal Processing}, vol.~66, no.~5, pp.
  1242--1257, March 2018.

\bibitem{ahlswede1986hypothesis}
R.~Ahlswede and I.~Csisz{\'a}r, ``Hypothesis testing with communication
  constraints,'' \emph{IEEE Transactions on Information Theory}, vol.~32,
  no.~4, pp. 533--542, 1986.

\bibitem{han1989exponential}
T.~S. Han and K.~Kobayashi, ``Exponential-type error probabilities for
  multiterminal hypothesis testing,'' \emph{IEEE Transactions on Information
  Theory}, vol.~35, no.~1, pp. 2--14, 1989.

\bibitem{watanabe2018neyman}
S.~Watanabe, ``Neyman--pearson test for zero-rate multiterminal hypothesis
  testing,'' \emph{IEEE Transactions on Information Theory}, vol.~64, no.~7,
  pp. 4923--4939, 2018.

\bibitem{10.5555/524893}
P.~K. Varshney, \emph{Distributed Detection and Data Fusion}, 1st~ed.\hskip 1em
  plus 0.5em minus 0.4em\relax Berlin, Heidelberg: Springer-Verlag, 1996.

\bibitem{4102537}
R.~R. {Tenney} and N.~R. {Sandell}, ``Detection with distributed sensors,''
  \emph{IEEE Transactions on Aerospace and Electronic Systems}, vol. AES-17,
  no.~4, pp. 501--510, 1981.

\bibitem{395226}
W.~{Baek} and S.~{Bommareddy}, ``Optimal m-ary data fusion with distributed
  sensors,'' \emph{IEEE Transactions on Aerospace and Electronic Systems},
  vol.~31, no.~3, pp. 1150--1152, 1995.

\bibitem{5722050}
C.~{Tepedelenlioglu} and S.~{Dasarathan}, ``Distributed detection over gaussian
  multiple access channels with constant modulus signaling,'' \emph{IEEE
  Transactions on Signal Processing}, vol.~59, no.~6, pp. 2875--2886, 2011.

\bibitem{4205085}
J.~{Chamberland} and V.~V. {Veeravalli}, ``Wireless sensors in distributed
  detection applications,'' \emph{IEEE Signal Processing Magazine}, vol.~24,
  no.~3, pp. 16--25, 2007.

\bibitem{4099562}
G.~{Mergen}, V.~{Naware}, and L.~{Tong}, ``Asymptotic detection performance of
  type-based multiple access over multiaccess fading channels,'' \emph{IEEE
  Transactions on Signal Processing}, vol.~55, no.~3, pp. 1081--1092, 2007.

\bibitem{7405352}
V.~S.~S. {Nadendla} and P.~K. {Varshney}, ``Design of binary quantizers for
  distributed detection under secrecy constraints,'' \emph{IEEE Transactions on
  Signal Processing}, vol.~64, no.~10, pp. 2636--2648, 2016.

\bibitem{698826}
R.~{Chandramouli} and N.~{Ranganathan}, ``Quantization for robust sequential
  m-ary signal detection,'' in \emph{ISCAS '98. Proceedings of the 1998 IEEE
  International Symposium on Circuits and Systems (Cat. No.98CH36187)}, vol.~4,
  1998, pp. 317--320 vol.4.

\bibitem{wigger2016testing}
M.~Wigger and R.~Timo, ``Testing against independence with multiple decision
  centers,'' in \emph{2016 International Conference on Signal Processing and
  Communications (SPCOM)}.\hskip 1em plus 0.5em minus 0.4em\relax IEEE, 2016,
  pp. 1--5.

\bibitem{xiang2012interactive}
Y.~Xiang and Y.-H. Kim, ``Interactive hypothesis testing with communication
  constraints,'' in \emph{2012 50th Annual Allerton Conference on
  Communication, Control, and Computing (Allerton)}.\hskip 1em plus 0.5em minus
  0.4em\relax IEEE, 2012, pp. 1065--1072.

\bibitem{salehkalaibar2018hypothesis}
S.~Salehkalaibar, M.~Wigger, and R.~Timo, ``On hypothesis testing against
  conditional independence with multiple decision centers,'' \emph{IEEE
  Transactions on Communications}, vol.~66, no.~6, pp. 2409--2420, 2018.

\bibitem{4558051}
K.~R. {Varshney} and L.~R. {Varshney}, ``Quantization of prior probabilities
  for hypothesis testing,'' \emph{IEEE Transactions on Signal Processing},
  vol.~56, no.~10, pp. 4553--4562, 2008.

\bibitem{6203607}
J.~B. {Rhim}, L.~R. {Varshney}, and V.~K. {Goyal}, ``Quantization of prior
  probabilities for collaborative distributed hypothesis testing,'' \emph{IEEE
  Transactions on Signal Processing}, vol.~60, no.~9, pp. 4537--4550, 2012.

\bibitem{nakagawa1993converse}
K.~Nakagawa and F.~Kanaya, ``On the converse theorem in statistical hypothesis
  testing,'' \emph{IEEE Transactions on Information Theory}, vol.~39, no.~2,
  pp. 623--628, 1993.

\bibitem{strassen2009asymptotic}
V.~Strassen, ``Asymptotic estimates in {S}hannon's information theory,'' in
  \emph{Proc. 3rd Trans. Prague Conf. Inf. Theory}, 2009, pp. 689--723.

\bibitem{sason2012moderate}
I.~Sason, ``Moderate deviations analysis of binary hypothesis testing,'' in
  \emph{2012 IEEE International Symposium on Information Theory
  Proceedings}.\hskip 1em plus 0.5em minus 0.4em\relax IEEE, 2012, pp.
  821--825.

\bibitem{dembo2011large}
A.~Dembo and O.~Zeitouni, ``Large deviations techniques and applications.
  1998,'' \emph{Applications of Mathematics}, vol.~38, 2011.

\bibitem{espinosa2021finite}
S.~Espinosa, J.~F. Silva, and P.~Piantanida, ``Finite-length bounds on
  hypothesis testing subject to vanishing type i error restrictions,''
  \emph{IEEE Signal Processing Letters}, no.~5, January 2021.

\bibitem{DBLP:journals/corr/KatzPD16a}
\BIBentryALTinterwordspacing
G.~Katz, P.~Piantanida, and M.~Debbah, ``Collaborative distributed hypothesis
  testing,'' \emph{CoRR}, vol. abs/1604.01292, 2016. [Online]. Available:
  \url{http://arxiv.org/abs/1604.01292}
\BIBentrySTDinterwordspacing

\bibitem{espinosa2019new}
S.~Espinosa, J.~F. Silva, and P.~Piantanida, ``New results on testing against
  independence with rate-limited constraints,'' in \emph{2019 IEEE Global
  Conference on Signal and Information Processing (GlobalSIP)}.\hskip 1em plus
  0.5em minus 0.4em\relax IEEE, 2019, pp. 1--5.

\bibitem{cover2012elements}
T.~M. Cover and J.~A. Thomas, \emph{Elements of information theory}.\hskip 1em
  plus 0.5em minus 0.4em\relax John Wiley \& Sons, 2012.

\bibitem{kendall1999kendall}
M.~Kendall, A.~Stuart, K.~J. Ord, and S.~Arnold, \emph{Kendall's Advanced
  Theory of Statistics: Volume 2A--Classical Inference and and the Linear
  Model}.\hskip 1em plus 0.5em minus 0.4em\relax Wiley Online Library, 1999.

\bibitem{kullback1951information}
S.~Kullback and R.~A. Leibler, ``On information and sufficiency,'' \emph{The
  Annals of Mathematical Statistics}, vol.~22, no.~1, pp. 79--86, 1951.

\bibitem{chernoff1952measure}
H.~Chernoff, ``A measure of asymptotic efficiency for tests of a hypothesis
  based on the sum of observations,'' \emph{The Annals of Mathematical
  Statistics}, pp. 493--507, 1952.

\bibitem{shkel}
Y.~Shkel, M.~Raginsky, and S.~Verd{\'u}, ``Universal lossy compression under
  logarithmic loss,'' in \emph{Information Theory (ISIT), 2017 IEEE
  International Symposium on}.\hskip 1em plus 0.5em minus 0.4em\relax IEEE,
  2017, pp. 1157--1161.

\bibitem{tishby2000information}
N.~Tishby, F.~C. Pereira, and W.~Bialek, ``The information bottleneck method,''
  \emph{arXiv preprint physics/0004057}, 2000.

\bibitem{han1987hypothesis}
T.~Han, ``Hypothesis testing with multiterminal data compression,'' \emph{IEEE
  Transactions on Information Theory}, vol.~33, no.~6, pp. 759--772, 1987.

\bibitem{boucheron2013concentration}
S.~Boucheron, G.~Lugosi, and P.~Massart, \emph{Concentration inequalities: A
  nonasymptotic theory of independence}.\hskip 1em plus 0.5em minus 0.4em\relax
  Oxford University Press, 2013.

\bibitem{zhang1997redundancy}
Z.~Zhang, E.-H. Yang, and V.~K. Wei, ``The redundancy of source coding with a
  fidelity criterion-part one: Known statistics,'' \emph{IEEE Transactions on
  Information Theory}, vol.~43, no.~1, pp. 71--91, 1997.

\bibitem{vera2018compression}
M.~Vera, L.~R. Vega, and P.~Piantanida, ``Compression-based regularization with
  an application to multitask learning,'' \emph{IEEE Journal of Selected Topics
  in Signal Processing}, vol.~12, no.~5, pp. 1063--1076, 2018.

\bibitem{arimoto1972algorithm}
S.~Arimoto, ``An algorithm for computing the capacity of arbitrary discrete
  memoryless channels,'' \emph{IEEE Transactions on Information Theory},
  vol.~18, no.~1, pp. 14--20, 1972.

\bibitem{gallager2008principles}
R.~Gallager, \emph{Principles of digital communication}.\hskip 1em plus 0.5em
  minus 0.4em\relax Technical Publications, 2008.

\bibitem{devroye2012combinatorial}
L.~Devroye and G.~Lugosi, \emph{Combinatorial methods in density
  estimation}.\hskip 1em plus 0.5em minus 0.4em\relax Springer Science \&
  Business Media, 2012.

\bibitem{courtade2014multiterminal}
T.~A. Courtade and T.~Weissman, ``Multiterminal source coding under logarithmic
  loss,'' \emph{IEEE Transactions on Information Theory}, vol.~60, no.~1, pp.
  740--761, 2014.

\bibitem{berger1971rate}
T.~Berger, ``Rate-distortion theory,'' \emph{Encyclopedia of
  Telecommunications}, 1971.

\bibitem{koga2013information}
T.~S. Han, \emph{Information-spectrum methods in information theory}.\hskip 1em
  plus 0.5em minus 0.4em\relax Springer Science \& Business Media, 2013.

\end{thebibliography}
\end{document}